\documentclass[reqno]{amsart}
\usepackage[T1]{fontenc}
\usepackage{amssymb}
\usepackage{amsthm}
\usepackage{thmtools}
\usepackage{xcolor}
\usepackage[colorlinks,allcolors=blue!50!black]{hyperref}

\newcommand{\norm}[2][y]{\if#1y\left\fi\lVert#2\if#1y\right\fi\rVert} 
\newcommand{\abs}[2][y]{\if#1y\left\fi\lvert#2\if#1y\right\fi\rvert} 
\renewcommand{\natural}{\mathbb{N}} 
\newcommand{\Real}{\mathbb{R}}
\newcommand{\Complex}{\mathbb{C}}
\newcommand{\LL}{\mathrm{L}}
\newcommand{\HH}{\mathrm{H}}
\newcommand{\CC}{\mathrm{C}}
\newcommand{\poly}{\Pi}
\newcommand{\dd}[1][y]{\if#1y\,\fi{\mathrm d}} 
\DeclareMathOperator{\proj}{proj}
\DeclareMathOperator{\vspan}{span}
\renewcommand{\div}{\operatorname{div}} 

\declaretheorem[name=Theorem,numberwithin=section]{thm}
\declaretheorem[name=Lemma,sibling=thm]{lemma}
\declaretheorem[name=Proposition,sibling=thm]{proposition}
\declaretheorem[name=Corollary,sibling=thm]{corollary}
\declaretheorem[name=Remark,sibling=thm,style=remark]{remark}

\begin{document}
\title[Sobolev orthogonal polynomials in the ball]{Weighted Sobolev orthogonal polynomials and approximation in the ball}
\author{Leonardo E. Figueroa}
\address{CI\textsuperscript{2}MA and Departamento de Ingenier\'ia Matem\'atica, Universidad de Concepci\'on, Concepci\'on, Chile}
\email{lfiguero@ing-mat.udec.cl}

\begin{abstract}
We establish simultaneous approximation properties of weighted first-order Sobolev orthogonal projectors onto spaces of polynomials of bounded total degree in the Euclidean unit ball.
The simultaneity is in the sense that we provide bounds for the projection error in both a weighted $\LL^2$ norm and a weighted $\HH^1$ seminorm, both involving the same weight of the generalized Gegenbauer type $x \mapsto (1-\norm{x}^2)^\alpha$, $\alpha > -1$.
The Sobolev orthogonal projectors producing the approximations are with respect to an alternative yet equivalent inner product for the corresponding uniformly weighted $\HH^1$ space.
In order to obtain our approximation bounds, we study the orthogonal polynomial structure of this alternative Sobolev inner product obtaining, among other results, a characterization of its orthogonal polynomials as solutions of certain Sturm--Liouville problems.
We do not rely on any particular basis of orthogonal polynomials, which allows for a streamlined and dimension-independent exposition.

\end{abstract}

\thanks{L.E.~Figueroa acknowledges funding from Vicerrector\'ia de Investigaci\'on y Desarrollo de la Universidad de Concepci\'on, project VRID Investigaci\'on \textnumero~220.013.047-INV}
\maketitle



\section{Introduction}

\subsection{Purpose}\label{ssc:purpose}
Let $\norm{\cdot}$ denote the Euclidean norm on $\Real^d$ and let $B^d$ be its unit ball.
Given $\alpha > -1$, let the weight function $W_\alpha \colon B^d \to \Real^d$ be defined by $W_\alpha(x) := (1-\norm{x}^2)^\alpha$.
Given an integer $N \geq 0$, let $\poly^d_N$ be the space of $d$-variate algebraic polynomials of total degree no higher than $N$.
Let also $\LL^2_\alpha$ be the weighted Lebesgue space $\LL^2_\alpha(B^d,W_\alpha) := \{ W_\alpha^{-1/2} f \mid f \in \LL^2(B^d) \}$, whose natural squared norm is $\norm{u}_\alpha^2 := \int_{B^d} \abs{u(x)}^2 W_\alpha(x) \dd x$.
Further, given an integer $m \geq 0$, let $\HH^m_\alpha$ be the weighted Sobolev space whose squared norm is $\norm{u}_{\alpha,m,\mathrm{S}}^2 := \sum_{k=0}^m \norm{\nabla^k u}_\alpha^2$ (the `S' stands for `standard'), with $\nabla^k$ denoting the weak $k$-fold gradient.

The purpose of this work is proving that, for each integer $m \geq 1$, each function $u \in \HH^{\alpha,m}$ and each degree $N$ that is high enough, there exists an algebraic polynomial $u_N \in \poly^d_N$ such that the following simultaneous approximation bounds hold:
\begin{subequations}\label{simultaneous}
\begin{gather}
\label{simultaneous-L2}%
\norm{u - u_N}_\alpha \leq C N^{-m} \norm{\nabla^m u}_\alpha,\\
\label{simultaneous-H1}%
\norm{\nabla u - \nabla u_N}_\alpha \leq C N^{1-m} \norm{\nabla^m u}_\alpha;
\end{gather}
\end{subequations}
here, $C$ is a positive constant that does not depend on the function $u$ nor the degree $N$ (see \autoref{cor:old-style-bounds} for the precise formulation of this result).

The approximants in \eqref{simultaneous} will be of the form $u_N = S^{\alpha,1}_N u$, where $S^{\alpha,1}_N \colon \HH^1_\alpha \to \poly^d_N$ is an orthogonal projector with respect to the alternative (yet equivalent in the sense of inducing the same Hilbert space topology) Sobolev inner product
\begin{equation}\label{projection-IP-1}
\langle u, v \rangle_{\alpha,1} := \langle \nabla u, \nabla v \rangle_\alpha + c_\alpha^{-1} \langle u, 1 \rangle_\alpha \overline{\langle v, 1 \rangle_\alpha},
\end{equation}
where, in turn, $\langle \cdot, \cdot \rangle_\alpha$ denotes the inner product of $\LL^2_\alpha$ and its vectorizations and tensorizations and $c_\alpha = \langle 1, 1 \rangle_\alpha$.
For this reason, the structure of orthogonal polynomials with respect to the inner product in \eqref{projection-IP-1} will be of paramount importance to us (which is why we give this non-standard inner product the shorter label).
Among its features that we will provide lie an orthogonal decomposition of the associated Sobolev orthogonal polynomial spaces (\autoref{lem:SOP}) and weak (\autoref{thm:SL-1}) and strong (\autoref{thm:second-order-SL-1}) Sturm--Liouville problems that are satisfied by their members.

We will employ these features to derive bounds on the projection error $u - S^{\alpha,1}_N u$ in the norm induced by the Sobolev inner product \eqref{projection-IP-1} (that is, a best approximation error; \autoref{thm:H1-approximation}) and in the $\LL^2_\alpha$ norm (\autoref{thm:L2-approximation}) in terms of other best approximation errors in both cases.
Then, we will recast those bounds in terms of seminorms as in \eqref{simultaneous} (\autoref{cor:old-style-bounds}).

We were originally drawn into the study of polynomial approximation results such as \eqref{simultaneous} for the weighted Sobolev space $\HH^1_\alpha$ because the latter is, in the $\alpha > 1$ case and up to simple isomorphisms, equal to certain Maxwellian-weighted Sobolev spaces on Euclidean balls naturally arising in the analysis of kinetic models of dilute polymers \cite[Eq.~(2.3b) \& Rem.~5.8]{FS:2012}.

\subsection{Related results}\label{ssc:related}
Approximants of the form $u_N = S^\alpha_N u$, where $S^\alpha_N \colon \LL^2_\alpha \to \poly^d_N$ is the standard $\LL^2_\alpha$ orthogonal projector, satisfy \eqref{simultaneous-L2} in isolation \cite[Cor.~2.4]{Figueroa:2017b}.
This can also be gleaned from \cite[Th.~1.1, Th.~1.2 \& Cor.~1.3]{PinarXu:2018}.
However, such $\LL^2_\alpha$ orthogonal projector behaves quantifiably worse than \eqref{simultaneous-H1} for the gradient; indeed, it follows from the case $r = 1$ of \cite[Th.~1.1]{Figueroa:2017a} that, for $N$ high enough,
\begin{equation*}
\norm{\nabla u - \nabla S^\alpha_N u}_\alpha \leq C \, N^{3/2-m} \norm{u}_{\alpha,m,\mathrm{S}}
\end{equation*}
for some positive constant $C$.
By standard arguments (cf.\ the proof of \autoref{cor:old-style-bounds}), one might replace the norm $\norm{u}_{\alpha,m,\mathrm{S}}$ with the $\norm{\nabla^m u}_\alpha$ seminorm in the above bound.
However, the power on $N$ cannot be lowered in the $m = 1$ case \cite[Th.~5.1]{Figueroa:2017a}, in which this bound actually encodes a divergence rate.

In \cite[Th.~4.1]{LiXu:2014} it is shown that, in the unweighted ($\alpha = 0$) case, approximants of the form $u_N = S^{-1}_N u$ satisfy the simultaneous approximation bounds \eqref{simultaneous}, where $S^{-1}_N \colon \HH^1(B^d) \to \poly^d_N$ is the orthogonal projection operator with respect to an equivalent inner product for the Sobolev space $\HH^1(B^d)$ proportional to
\begin{equation}\label{LiXu-IP}
(u,v) \mapsto \int_{B^d} \nabla u(x) \cdot \overline{\nabla v(x)} \, \dd x + \lambda_0 \int_{\mathbb{S}^{d-1}} u(x) \overline{v(x)} \dd\sigma_{d-1}(x),
\end{equation}
where $\lambda_0$ is an arbitrary positive constant, $\mathbb{S}^{d-1} = \partial B^d$ is the unit sphere and $\sigma_{d-1}$ is its standard surface measure.
The result \cite[Th.~4.1]{LiXu:2014} also provides bounds in terms of other best approximation errors.
This important work also encompasses approximation in Sobolev spaces $W^{m,p}(B^d)$ of higher regularity $m$ and general integrability parameter $p \in (1,\infty)$ \cite[Th.~4.1 and Th.~4.2]{LiXu:2014} and relies heavily on explicit bases of orthogonal polynomials with respect to \eqref{LiXu-IP} and other Sobolev inner products.

Analogous results concerning unweighted $\LL^2$ simultaneous approximation of a function and its first and second derivatives in the unit triangle, whose orthogonal structure is much more complicated than that of the ball, are obtained in \cite[Th.~1.2 \& Th.~1.3]{Xu:2017b}.

There is also \cite{MPPR:2023}, which provides error bounds for the best approximation in $\poly^d_n$ with respect to a Sobolev inner product proportional to $(u,v) \mapsto \rho \int_{B^d} \nabla u(x) \cdot \nabla v(x) W_\alpha(x) \dd x + u(0) v(0)$, $\rho > 0$ arbitrary and $\alpha \geq 0$.
As part of their argument, explicit bases of orthogonal polynomials with respect to that inner product are constructed in terms, ultimately, of Jacobi polynomials and any $\LL^2(\mathbb{S}^{d-1})$-orthogonal basis of spherical harmonics.

We now mention other works delving on the structure of multivariate Sobolev polynomials which, unlike the already mentioned \cite{LiXu:2014} and \cite{Xu:2017b}, are not explicitly concerned with approximation results.
We start with \cite{Xu:2006a}, where an orthogonal basis of Sobolev orthogonal polynomials in the ball with respect to an inner product proportional to $(u,v) \mapsto \int_{B^d} \Delta\mathopen{}\left( (1-\norm{x}^2) u(x) \mathclose{}\right) \Delta\mathopen{}\left( (1-\norm{x}^2) v(x) \mathclose{}\right) \dd[n] x$ is constructed in terms of Jacobi polynomials and any $\LL^2(\mathbb{S}^{d-1})$-orthogonal basis of spherical harmonics.

In \cite{Xu:2008a} similar bases of Sobolev orthogonal polynomials are constructed for inner products proportional to \eqref{LiXu-IP} and $(u,v) \mapsto \rho \int_{B^d} \nabla u(x) \cdot \nabla v(x) \dd x + u(0) v(0)$, where in both cases $\rho$ is an arbitrary positive constant.
This work also notes in \cite[Cor.~2.4]{Xu:2008a} that the space of orthogonal polynomials of degree $n \geq 1$ with respect to the inner product \eqref{LiXu-IP} admits an orthogonal decomposition of the form
\begin{equation}\label{Xu-OD}
\mathcal{H}^d_n \oplus (1-\norm{x}^2) \, \mathcal{V}^1_{n-2},
\end{equation}
where $\mathcal{H}^d_n$ is the space of $d$-variate spherical harmonics of degree $n$ and $\mathcal{V}^1_{n-2}$ is the space of $\LL^2_1$-orthogonal polynomials of degree $n-2$ (cf.\ \eqref{V0});
it further notes in \cite[Rem.~2.1]{Xu:2008a} the correct form of a similar decomposition for the spaces of orthogonal polynomials with respect to the inner product studied in \cite{Xu:2006a}.
This work goes on to provide in \cite[Cor.~2.5]{Xu:2008a} a fourth-order differential equation that is satisfied by the Sobolev orthogonal spaces with respect to the inner product \eqref{LiXu-IP}.

In \cite{PinarXu:2009} Sobolev orthogonal polynomials with respect to the inner product \eqref{LiXu-IP} are shown to satisfy the same second-order Sturm--Liouville problem $\mathcal{L}\sp{(\alpha)} y = {(\alpha)}_n y$ (cf.\ \eqref{SL-0-strong} and \eqref{SL-0-strong-ew}) that $W_\alpha$-Lebesgue orthogonal polynomials satisfy for $\alpha > -1$, but with $\alpha = -1$.
Families of polynomial solutions for the same Sturm--Liouville problem are also constructed for $\alpha = -2, -3, \dotsc$; in the $\alpha = -2$ case, an explicit perturbation of the corresponding family is shown to consist of Sobolev orthogonal polynomials with respect to the inner product $(u,v) \mapsto \rho \int_{B^d} \Delta u(x) \Delta v(x) \dd x + \int_{\mathbb{S}^{d-1}} u(x) v(x) \dd\sigma_{d-1}(x)$, $\rho > 0$ arbitrary \cite[Th.~4.1]{PinarXu:2009}.

In \cite{PPX:2013} orthogonal bases of Sobolev orthogonal polynomials with respect to a standard weighted Sobolev inner product proportional to $(u,v) \mapsto \rho \int_{B^d} \nabla u(x) \cdot \nabla v(x) W_\alpha(x) \dd x + \int_{B^d} u(x) v(x) W_\alpha(x) \dd x$, $\rho > 0$ arbitrary, are constructed for all $\alpha > -1$ in terms of any $\LL^2(\mathbb{S}^{d-1})$-orthogonal basis of spherical harmonics and new families of univariate Sobolev orthogonal polynomials.
There is also a brief aside to construct orthogonal bases of Sobolev orthogonal polynomials with respect to inner products proportional to $(u,v) \mapsto \rho \int_{B^d} \nabla u(x) \cdot \nabla v(x) W_{\alpha+1}(x) \dd x + \int_{\mathbb{S}^{d-1}} u(x) v(x) \dd \sigma_{d-1}(x)$, $\rho > 0$ arbitrary and $\alpha > -1$, in terms of Jacobi polynomials and any $\LL^2(\mathbb{S}^{d-1})$-orthogonal basis of spherical harmonics.

In \cite{DFLPP:2016}, orthogonal bases of Sobolev orthogonal polynomials in the ball with respect to inner products proportional to $(u,v) \mapsto \int_{B^d} u(x) v(x) W_\alpha(x) \dd x + \rho \int_{\mathbb{S}^{d-1}} \frac{\partial u}{\partial \nu}(x) \frac{\partial v}{\partial \nu}(x) \dd\sigma_{d-1}(x)$, $\rho > 0$ arbitrary, $\alpha > -1$, where $\frac{\partial}{\partial\nu} = x \cdot \nabla$ is the outward normal derivative operator, are constructed in terms of univariate Sobolev orthogonal polynomials with respect to a perturbation of the standard Jacobi inner product and any $\LL^2(\mathbb{S}^{d-1})$-orthogonal basis of spherical harmonics.

In \cite{LPP:2021}, orthogonal bases of Sobolev orthogonal polynomials in the ball with respect to inner products proportional to $(u,v) \mapsto \int_{B^d} u(x) v(x) W_\alpha(x) \dd x + \rho \int_{B^d} (x \cdot \nabla u(x)) (x \cdot \nabla v(x)) W_\alpha(x) \dd x$, $\rho > 0$ arbitrary, $\alpha > -1$, are constructed in terms of univariate Sobolev orthogonal polynomials with respect to another perturbation of the standard Jacobi inner product and any $\LL^2(\mathbb{S}^{d-1})$-orthogonal basis of spherical harmonics.

\subsection{Additional preliminary definitions and outline of the rest of this work}

We will denote the set of strictly positive integers by $\natural$ and let $\natural_0 := \{0\} \cup \natural$.
We denote the canonical basis of $\Real^d$ by $\{e_1, \dotsc, e_d\}$; that is, for $1 \leq i, j \leq d$, the $j$-th component of $e_i$ is $1$ if $i = j$ and $0$ otherwise.
Given a multi-index $\gamma \in [\natural_0]^d$ and a sufficiently (strongly or weakly) differentiable function $f$ we shall write $\abs{\gamma} = \sum_{i=1}^d \gamma_i$ and $\partial_\gamma f = \partial^{\abs{\gamma}}/(\partial x_1^{\gamma^1} \dotsm \partial x_d^{\gamma_d})$; we also use $\partial_i := \partial_{e_i} = \partial/\partial x_i$, $1 \leq i \leq d$.
We adopt the convention $\poly^d_n = \{0\}$ for $n < 0$.

Given $m \in \natural_0$, the inner product $\langle u, v \rangle_{\alpha,m,\mathrm{S}}$ on $\HH^m_\alpha$ defined by
\begin{equation}\label{typical-IP}
\langle u, v \rangle_{\alpha,m,\mathrm{S}}
:= \sum_{k=0}^m \langle \nabla^k u, \nabla^k v \rangle_\alpha
\end{equation}
is of course the inner product that induces the standard weighted Sobolev inner product $\norm{\cdot}_{\alpha,m,\mathrm{S}}$.
This inner product is equivalent to the commonly used $(u,v) \mapsto \sum_{\abs{\gamma} \leq m} \langle \partial_\gamma u, \partial_\gamma v \rangle_\alpha$.

The rest of this work is organized as follows.
In \autoref{sec:OP} we mostly collect known properties of Lebesgue orthogonal polynomial spaces with respect to the $\LL^2_\alpha$-inner product, including characterizations as eigenspaces of Sturm--Liouville operators and prove a variant of a known approximation result.
\autoref{sec:SOP} and \autoref{sec:approximation} form the core of this work.
In \autoref{sec:SOP} we introduce Sobolev orthogonal polynomial spaces with respect to the $(\alpha,1)$-inner product of \eqref{projection-IP-1}, derive a useful orthogonal decomposition and characterize them as eigenspaces of new Sturm--Liouville operators.
In \autoref{sec:approximation} we prove bounds for the error committed by the corresponding Sobolev orthogonal projection operators onto polynomials in terms of other best approximation errors and weave them into the desired simultaneous approximation bounds \eqref{simultaneous}.
In \autoref{sec:miscellanea} we succinctly discuss certain aspects of the orthogonal structure of polynomials with respect to the $(\alpha,1)$-inner product that do not play a direct role in the core of this work, but feature repeatedly in the literature cited above in \autoref{ssc:related}; then, we finish with a brief conclusion.

\section{Lebesgue orthogonal polynomials}\label{sec:OP}

Given $n \in \natural_0$, we denote by $\mathcal{V}^\alpha_n$ the space of $\LL^2_\alpha$-orthogonal polynomials of degree $n$; that is,
\begin{equation}\label{V0}
\mathcal{V}^\alpha_n := \{ p \in \poly^d_n \mid (\forall\,q\in \poly^d_{n-1})\ \langle p, q \rangle_\alpha = 0 \}.
\end{equation}
We adopt the convention $\mathcal{V}^\alpha_n = \{0\}$ for $n < 0$.
There holds (cf.\ \cite[Sec.~3.1]{DunklXu:2014})
\begin{equation}\label{dim-OP}
(\forall\,n\in\natural_0) \quad \dim(\mathcal{V}^\alpha_n) = \binom{n+d-1}{n}.
\end{equation}
For every $\alpha > -1$, as $W_\alpha$ is centrally symmetric, from \cite[Th.~3.3.11]{DunklXu:2014} there follows the parity relation
\begin{equation}\label{parity}
(\forall\,n\in\natural_0)\ (\forall\,p_n \in \mathcal{V}^\alpha_n)\ (\forall\,x\in B^d) \quad p_n(-x) = (-1)^n \, p_n(x).
\end{equation}
We denote by $\proj^\alpha_n$ the orthogonal projection from $\LL^2_\alpha$ onto $\mathcal{V}^\alpha_n$ and we have already introduced in \autoref{ssc:related} the notation $S^\alpha_n$ for the orthogonal projection from $\LL^2_\alpha$ onto $\poly^d_n$.
We will denote the entrywise application of $S^\alpha_n$ to members of vector- or higher-order tensor-valued variants of $\LL^2_\alpha$ by $S^\alpha_n$ as well.
From \cite[Th.~3.2.18]{DunklXu:2014}, $\poly^d_n = \bigoplus_{k=0}^n \mathcal{V}^\alpha_k$ and $\LL^2_\alpha = \bigoplus_{k=0}^\infty \mathcal{V}^\alpha_k$, whence
\begin{equation*}
(\forall\,n\in\natural_0)\quad S^\alpha_n = \sum_{k=0}^n \proj^\alpha_k
\qquad\text{and}\qquad
(\forall\,u\in\LL^2_\alpha)\quad u = \sum_{k=0}^\infty \proj^\alpha_k(u),
\end{equation*}
the latter series converging in $\LL^2_\alpha$.

Let us introduce the first-order angular differential operator \cite[Sec.~1.8]{DaiXu:2013}
\begin{equation}\label{Dij}
D_{i,j} := x_i \partial_j - x_j \partial_i
\end{equation}
Note that $D_{j,i} = -D_{i,j}$; in particular, $D_{i,i} = 0$.
It is immediate to check that these operators satisfy the relations $D_{i,j}(f\,g) = D_{i,j} f \, g + f \, D_{i,j} g$ and that they vanish on radial functions.

Given any $\alpha \in \Real$ and $j \in \{1, \dotsc, d\}$ let the first-order differentiation operator $d^\alpha_j$ be defined by
\begin{multline}\label{PLO}
d^\alpha_j q(x) := -(1-\norm{x}^2)^{-\alpha} \frac{\partial}{\partial x_j} \left( (1-\norm{x}^2)^{\alpha+1} q(x) \right)\\
= -(1-\norm{x}^2) \, \partial_j q(x) + 2 (\alpha+1) \, x_j \, q(x).
\end{multline}
Clearly, for all $\alpha \in \Real$ and $j \in \{1, \dotsc, d\}$, $d^\alpha_j \poly^d_n \subseteq \poly^d_{n+1}$.

Given $\alpha \in \Real$, we define the second-order differential operator $\mathcal{L}\sp{(\alpha)}$ by
\begin{multline}\label{L}
\mathcal{L}\sp{(\alpha)}(q) := -W_\alpha^{-1} \div\left( W_{\alpha+1} \nabla q \right) - \sum_{1 \leq i<j \leq d} (x_i\partial_j - x_j\partial_i)^2 q\\
\stackrel{\text{\eqref{Dij},\eqref{PLO}}}{=} \sum_{j=1}^d d^\alpha_j \partial_j q - \sum_{1 \leq i<j \leq d} D_{i,j}^2 q.
\end{multline}
Members of $\mathcal{V}^\alpha_n$, $\alpha>-1$, satisfy the second-order Sturm--Liouville problem \cite[Eq.~(5.2.3) and Th.~8.1.3]{DunklXu:2014}
\begin{equation}\label{SL-0-strong}
(\forall\,p_n \in \mathcal{V}^\alpha_n) \quad \mathcal{L}\sp{(\alpha)}(p_n) = \lambda\sp{(\alpha)}_n \, p_n,
\end{equation}
where
\begin{equation}\label{SL-0-strong-ew}
\lambda\sp{(\alpha)}_n = n(n+d+2\alpha).
\end{equation}

The angular differential operator $D_{i,j}$ is minus its own adjoint in the sense of the following result (see \cite[Prop.~1.8.4]{DaiXu:2013} for a variant on the unit sphere).

\begin{proposition}\label{pro:Dij-IBP} Let $d \in \natural$, $\alpha > -1$, $f, g \in \CC^1(\overline{B^d})$.
Then,
\begin{equation*}
\langle D_{i,j} f, g \rangle_\alpha = - \langle f, D_{i,j} g \rangle_\alpha.
\end{equation*}
\begin{proof}
As the weight $W_\alpha$ is a radial function,
\begin{equation*}
\left( D_{i,j} f \, \overline{g} + f \, \overline{D_{i,j} g} \right) W_\alpha
= \div\mathopen{}\left( f \, \overline{g} \, W_\alpha \, (x_i e_j - x_j e_i) \right)\mathclose{}
\end{equation*}
in the open unit ball $B^d$.
Let $r \in (0,1)$.
Because of Lebesgue's dominated convergence theorem, the integral of the left-hand side above over $B^d$ is the one-sided limit as $r \to 1^-$ of the corresponding integral over $r B^d$, which, on account of the divergence theorem, vanishes.
\end{proof}
\end{proposition}

Given any polynomial $p$ and $i,j \in \{1, \dotsc, d\}$, clearly $D_{i,j} p$ can only be a polynomial of the same degree as $p$ or null.
Combining this fact with \autoref{pro:Dij-IBP} we find that, for all $\alpha > -1$, there holds the inclusion
\begin{equation}\label{Dij-no-shift}
D_{i,j} \mathcal{V}^\alpha_n \subseteq \mathcal{V}^\alpha_n.
\end{equation}

The other two operators comprising the Sturm--Liouville operator $\mathcal{L}\sp{(\alpha)}$ of \eqref{L}, when $\alpha > -1$, satisfy another adjointness relation \cite[Prop.~3.2(ii)]{Figueroa:2017a}:
For $f, g \in \CC^1(\overline{B^d})$,
\begin{equation}\label{PLO-adjoint}
\langle \partial_j f, g \rangle_{\alpha+1} = \langle f, d^\alpha_j g \rangle_\alpha.
\end{equation}
Still under the assumption that $\alpha > -1$, these two operators also satisfy the inclusions \cite[Prop.~3.2(iii)]{Figueroa:2017a}
\begin{equation}\label{parameter-lowering}
d^\alpha_j \mathcal{V}^{\alpha+1}_n \subseteq \mathcal{V}^\alpha_{n+1}.
\end{equation}
and \cite[Lem.~2.11]{LiXu:2014}
\begin{equation}\label{diff-shift}
\partial_j \mathcal{V}^\alpha_n \subseteq \mathcal{V}^{\alpha+1}_{n-1}.
\end{equation}

With the help of \autoref{pro:Dij-IBP} and \eqref{PLO-adjoint}, for $\alpha > -1$ the Sturm--Liouville operator of \eqref{L} can be shown to be self-adjoint with respect to $\LL^2_\alpha$; indeed, for all $p \in \CC^2(\overline{B^d})$ and $q \in \CC^1(\overline{B^d})$,
\begin{equation}\label{L-self-adjoint}
\langle \mathcal{L}\sp{(\alpha)}p, q \rangle_\alpha = \langle \nabla p, \nabla q \rangle_{\alpha+1} + \sum_{1\leq i<j\leq d} \langle D_{i,j}p, D_{i,j} q \rangle_\alpha =: B^\alpha(p,q).
\end{equation}
Thus, the Sturm--Liouville problem \eqref{SL-0-strong} satisfied by $\LL^2_\alpha$-orthogonal polynomials can be immediately recast into the weak form (cf.\ \cite[Eq.~(6)]{Figueroa:2017a}):
\begin{equation}\label{SL-0-weak}
(\forall\,p_n \in \mathcal{V}^{\alpha}_n)\ (\forall\,q\in\CC^1(\overline{B^d}))\quad
B^\alpha(p_n,q)
= \lambda\sp{(\alpha)}_n \, \langle p_n, q \rangle_\alpha.
\end{equation}
Let us note that the sesquilinear form $B^\alpha$ defined in \eqref{L-self-adjoint} still makes sense if its arguments lie in $\HH^1_\alpha$.
It is straightforward to check that, for every $u \in \HH^1_\alpha$ and almost every $x \in B^d$,
\begin{equation*}
(1-\norm{x}^2) \norm{\nabla u(x)}^2 + \sum_{1\leq i<j\leq d} \abs{D_{i,j} u(x)}^2
= \norm{\nabla u(x)}^2 - \abs{x\cdot\nabla u(x)}^2
\leq \norm{\nabla u(x)}^2,
\end{equation*}
so
\begin{equation}\label{B-bound}
(\forall\,u,v\in\HH^1_\alpha) \quad \abs{B^\alpha(u,v)} \leq \norm{\nabla u}_\alpha \norm{\nabla v}_\alpha.
\end{equation}

We recall that $\mathcal{H}^d_n$ denotes the space of $d$-variate spherical harmonics (i.e., homogeneous and harmonic polynomials) of degree $n$ \cite[Sec.~1.1]{DaiXu:2013}, \cite[Sec.~4.1]{DunklXu:2014}.
The members of the $\mathcal{H}^d_n$ are in a one-to-one relation with their restrictions to the unit sphere $\mathbb{S}^{d-1}$ and to the latter is the term `spherical harmonic' commonly reserved; we won't follow that practice.
We adopt the convention $\mathcal{H}^d_n = \{0\}$ if $n < 0$.
From \cite[Th.~4.1.2]{DunklXu:2014},
\begin{equation}\label{dim-SH}
\dim\left(\mathcal{H}^d_n\right) = \binom{n+d-1}{d-1} - \binom{n+d-3}{d-1}.
\end{equation}
For $n \in \natural_0$, let $\{Y^n_\nu\}_{\nu=1}^{\dim(\mathcal{H}^d_n)}$ be an $\LL^2(\mathbb{S}^{d-1})$-orthogonal basis of $\mathcal{H}^d_n$.
Also, let $P^{(\alpha,\beta)}_n$ denote the Jacobi polynomial of parameter $(\alpha,\beta)$ and degree $n$ \cite[Subsec.~1.4.4]{DunklXu:2014}.

\begin{proposition}[{\cite[Prop.~5.2.1]{DunklXu:2014}}]\label{pro:OP-basis}
Let $d \in \natural$, $\alpha > -1$ and $n \in \natural_0$.
Then, the polynomials $P^n_{j,\nu}$, where $j \in \{0, \dotsc, \lfloor n/2 \rfloor\}$ and $\nu \in \{1, \dotsc, \dim(\mathcal{H}^d_{n-2j})\}$, defined by
\begin{equation*}
P^n_{j,\nu}(x) := P^{\left(\alpha,n-2j+\frac{d-2}{2}\right)}_j(2\norm{x}^2-1) \, Y^{n-2j}_\nu(x),
\end{equation*}
form an $\LL^2_\alpha(B^d)$-orthogonal basis of $\mathcal{V}^\alpha_n$.
\end{proposition}

As every Jacobi polynomial of degree $0$ is the constant function $1$, it follows from \autoref{pro:OP-basis} that, for all $\alpha > -1$,
\begin{equation}\label{SH-are-OP}
(\forall\,n\in\natural_0) \quad \mathcal{H}^d_n \subseteq \mathcal{V}^\alpha_n
\end{equation}
with equality if $n = 0$ or $n = 1$, so
\begin{equation}\label{low-degree-OP}
\mathcal{V}^\alpha_0 = \vspan(\{1\})
\quad\text{and}\quad
\mathcal{V}^\alpha_1 = \vspan\left((x \mapsto x_i)_{i=1}^d\right).
\end{equation}

We finish this section with a variant of Theorem~1.1 and Theorem~1.2 of \cite{PinarXu:2018} where the $\LL^2_\alpha$-best approximation error of a regular enough function is bounded in terms of best approximation errors of its gradient and of its image under the Sturm--Liouville operator $\mathcal{L}\sp{(\alpha)}$.
First, we state a density result.

\begin{proposition}[{\cite[Lemma~2.2]{Figueroa:2017a}}]\label{pro:C-infty-density}
Let $d \in \natural$ and $\alpha > -1$.
Then, $\CC^\infty(\overline{B^d})$ is dense in $\HH^m_\alpha$.
\end{proposition}

\begin{thm}\label{thm:proj-L2-res-L2}
Let $d \in \natural$ and $\alpha > -1$.
Then,
\begin{equation}\label{proj-L2-grad-BA}
(\forall\,u\in\HH^1_\alpha)\ (\forall\,N\in\natural_0) \quad
\norm{u - S^\alpha_N(u)}_\alpha \leq (\lambda\sp{(\alpha)}_{N+1})^{-1/2} \inf_{p_{N} \in \poly^d_{N}} \norm{\nabla u - \nabla p_N}_\alpha
\end{equation}
and
\begin{equation}\label{proj-L2-SL-0-BA}
(\forall\,u\in\HH^2_\alpha)\ (\forall\,N\in\natural_0) \quad
\norm{u - S^\alpha_N(u)}_\alpha \leq (\lambda\sp{(\alpha)}_{N+1})^{-1} \norm{\mathcal{L}\sp{(\alpha)}(u) - S^\alpha_N(\mathcal{L}\sp{(\alpha)}(u))}_\alpha.
\end{equation}
\begin{proof}
Let $v \in \HH^1_\alpha$.
By \autoref{pro:C-infty-density}, it is the limit in $\HH^1_\alpha$ of a sequence of $\CC^1(\overline{B^d})$ functions.
Hence, exploiting the bound \eqref{B-bound} and the structure of the Sobolev inner product \eqref{typical-IP}, we can extend \eqref{SL-0-weak} to:
\begin{equation*}
(\forall\,p_n\in\mathcal{V}^\alpha_n) \quad
B^\alpha(p_n,v) = \lambda\sp{(\alpha)}_n \langle p_n, v \rangle_\alpha.
\end{equation*}
Then, given any $N \in \natural_0$ and any polynomial $q \in \poly^d_N$,
\begin{equation*}
B^\alpha(q,v) = \sum_{n=0}^N B^\alpha(\proj^\alpha_n(q),v)
= \sum_{n=0}^N \lambda\sp{(\alpha)}_n \langle \proj^\alpha_n(q), v \rangle_\alpha,
\end{equation*}
so
\begin{equation}\label{B-diff-diff}
B^\alpha(v-q,v-q) = B^\alpha(v,v) + \sum_{n=0}^N \lambda\sp{(\alpha)}_n \left[ \norm{\proj^\alpha_n(q)}_\alpha^2 - 2 \Re \langle \proj^\alpha_n(q), v \rangle_\alpha \right].
\end{equation}
As the eigenvalues $\lambda\sp{(\alpha)}_n$ are nonnegative \eqref{SL-0-strong-ew}, a minimizer of the left-hand side of \eqref{B-diff-diff} among all $q \in \poly^d_N$ is obtained by minimizing what lies inside the square brackets of the right-hand side for each $n$ independently, and that is attained by choosing $q$ so that $\proj^\alpha_n(q) = \proj^\alpha_n(v)$ for $0 \leq n \leq N$; hence,
\begin{equation}\label{minimizer}
(\forall\,v\in\HH^1_\alpha)\ \left(\forall\,N\in\natural_0\right) \quad S^\alpha_N(v) \in \operatorname*{arg\,min}_{q\in\poly^d_N} B^\alpha(v-q,v-q).
\end{equation}
Now, for all $N \in \natural_0$,
\begin{equation*}
0 \leq B^\alpha(v - S^\alpha_N(v), v - S^\alpha_N(v))
\stackrel{\text{\eqref{B-diff-diff}}}{=} B^\alpha(v,v) - \sum_{n=0}^N \lambda\sp{(\alpha)}_n \norm{\proj^\alpha_n(v)}_\alpha^2.
\end{equation*}
Thus, every partial sum of the series of non-negative terms $\sum_{n=0}^\infty \lambda\sp{(\alpha)}_n \norm{\proj^\alpha_n(v)}_\alpha^2$ is bounded by the finite quantity $B^\alpha(v,v)$, so the series converges and we obtain the Bessel-type bound
\begin{equation}\label{Bessel-B}
(\forall\,v\in\HH^1_\alpha)\ (\forall\,N\in\natural_0) \quad \sum_{n=0}^\infty \lambda\sp{(\alpha)}_n \norm{\proj^\alpha_n(v)}_\alpha^2 \leq B^\alpha(v,v).
\end{equation}
Therefore, given $u \in \HH^1_\alpha$ and $N \in \natural_0$,
\begin{multline*}
\norm{u - S^\alpha_N(u)}_\alpha^2
= \sum_{n=N+1}^\infty \norm{\proj^\alpha_n(u)}_\alpha^2
\leq \sup_{n\geq N+1} (\lambda\sp{(\alpha)}_n)^{-1} \sum_{n=N+1}^\infty \lambda\sp{(\alpha)}_n \norm{\proj^\alpha_n(u)}_\alpha^2\\
\stackrel{\text{\eqref{Bessel-B}}}{\leq} (\lambda\sp{(\alpha)}_{N+1})^{-1} B^\alpha(u-S^\alpha_N(u), u-S^\alpha_N(u))\\
\stackrel{\text{\eqref{minimizer}}}{=} (\lambda\sp{(\alpha)}_{N+1})^{-1} \min_{p_N \in \poly^d_N} B^\alpha(u-p_N,u-p_N)
\stackrel{\text{\eqref{B-bound}}}{\leq} (\lambda\sp{(\alpha)}_{N+1})^{-1} \inf_{p_N \in \poly^d_N} \norm{\nabla u - \nabla p_N}_\alpha^2,
\end{multline*}
whence \eqref{proj-L2-grad-BA}.

Let $w \in \CC^2(\overline{B^d})$.
From \eqref{SL-0-strong}, \eqref{SL-0-strong-ew}, \eqref{L-self-adjoint} and the self-adjointness of $B^\alpha$ made evident there, for all $n \in \natural_0$ and $q_n \in \mathcal{V}^\alpha_n$,
\begin{equation*}
\langle \mathcal{L}\sp{(\alpha)} w, q_n \rangle_\alpha
= \lambda\sp{(\alpha)}_n \langle w, q_n \rangle_\alpha
= \langle \lambda\sp{(\alpha)}_n \proj^\alpha_n(w), q_n \rangle_\alpha.
\end{equation*}
So, let $u \in \HH^2_\alpha$ and $N \in \natural_0$.
The above equality, the fact that $\mathcal{L}\sp{(\alpha)}$ is a continuous map from $\HH^2_\alpha$ to $\LL^2_\alpha$ (cf.\ the second form in \eqref{L}) and the density result in \autoref{pro:C-infty-density} give
\begin{equation}\label{L-proj}
(\forall\,n\in\natural_0) \quad \proj^\alpha_n(\mathcal{L}\sp{(\alpha)} u) = \lambda\sp{(\alpha)}_n \, \proj^\alpha_n(u).
\end{equation}
Then, the bound \eqref{proj-L2-SL-0-BA} follows from
\begin{multline*}
\norm{u - S^\alpha_N(u)}_\alpha^2 \leq \sup_{n \geq N+1} (\lambda\sp{(\alpha)}_n)^{-2} \sum_{n=N+1}^\infty \norm{\lambda\sp{(\alpha)}_n \, \proj^\alpha_n(u)}_\alpha^2\\
\stackrel{\text{\eqref{L-proj}}}{=} (\lambda\sp{(\alpha)}_{N+1})^{-2} \sum_{n=N+1}^\infty \norm[n]{\proj^\alpha_n(\mathcal{L}\sp{(\alpha)} u)}_\alpha^2
= (\lambda\sp{(\alpha)}_{N+1})^{-2} \norm[n]{\mathcal{L}\sp{(\alpha)} u - S^\alpha_N(\mathcal{L}\sp{(\alpha)} u)}_\alpha^2.
\end{multline*}
\end{proof}
\end{thm}

\section{Sobolev orthogonal polynomials}\label{sec:SOP}

In \autoref{thm:proj-L2-res-L2} we exploited the strong \eqref{SL-0-strong} and weak \eqref{SL-0-weak} forms of the Sturm--Liouville problem that the $\LL^2_\alpha$-orthogonal polynomials satisfy to obtain the bounds \eqref{proj-L2-grad-BA} and \eqref{proj-L2-SL-0-BA}, which in turn can be used to prove a bound of the type \eqref{simultaneous-L2} (cf.\ \cite[Cor.~2.4]{Figueroa:2017b}).
The purpose of this \autoref{sec:SOP} is proving that polynomials orthogonal with respect to the $(\alpha,1)$-inner product for $\HH^1_\alpha$ introduced in \eqref{projection-IP-1} satisfy their own strong and weak Sturm--Liouville problems, so that later, in \autoref{sec:approximation}, we can reproduce the argument of \autoref{thm:proj-L2-res-L2} to obtain analogous bounds for an $\HH^1_\alpha$-best approximation error that, in turn, will be used to obtain the simultaneous approximation bounds \eqref{simultaneous} in full.

Recalling that $S^\alpha_n$ is the $\LL^2_\alpha$-orthogonal projection onto $\poly^d_n$, that $\mathcal{V}^\alpha_0 = \vspan(\{1\})$ and $c_\alpha = \langle 1, 1 \rangle_\alpha$, we can express the $(\alpha,1)$-inner product as
\begin{equation}\label{main-IP}
\langle u, v \rangle_{\alpha,1} = \langle \nabla u, \nabla v \rangle_\alpha + \langle S^\alpha_0(u), S^\alpha_0(v) \rangle_\alpha,
\end{equation}
which will be our preferred form.
Roughly speaking, the lower-order term of \eqref{main-IP} will mostly get out of the way; by way of contrast, in \cite{PPX:2013}, concerning the construction of bases of orthogonal polynomials for the first-order standard Sobolev inner product $\langle \nabla u, \nabla v \rangle_\alpha + \langle u, v \rangle_\alpha$, it is stated that
\begin{quote}
Likely, it is this tangle between the two terms that makes this case so much more complicated than [the case of \eqref{LiXu-IP}].
\end{quote}
We do not introduce a term-balancing positive constant before any of the two terms of \eqref{main-IP} because, as will be proved later (\autoref{ssc:balancing-constant}), doing so does not alter its orthogonal polynomial spaces nor the orthogonal projection operators onto the latter.
This inner product \eqref{main-IP} belongs to the wider family \eqref{projection-IP} of inner products, all of whose members will prove useful later.

\begin{proposition}\label{pro:projection-IP}
Let $d \in \natural$, $\alpha > -1$ and $m \in \natural_0$.
Then,
\begin{equation}\label{projection-IP}
\langle u, v \rangle_{\alpha,m} := \langle \nabla^m u, \nabla^m v \rangle_\alpha + \langle S^\alpha_{m-1}(u), S^\alpha_{m-1}(v) \rangle_\alpha
\end{equation}
defines an inner product for $\HH^m_\alpha$ equivalent to the standard $(\alpha,1,\mathrm{S})$-inner product given in \eqref{typical-IP}.
\begin{proof}
This is \cite[Prop.~2.6]{Figueroa:2017b} with the help of the fact that $\langle \nabla^m u, \nabla^m v \rangle_\alpha$ is equivalent to $\sum_{\abs{\gamma}=m} \langle \partial_\gamma u, \partial_\gamma v \rangle_\alpha$.
\end{proof}
\end{proposition}

\begin{remark}\label{rem:CompletenessAndTraces}
We remark in passing that the $\HH^\alpha_m$ spaces, with their intrinsic definition as given in \autoref{ssc:purpose}, are indeed Hilbert spaces \cite{KO}.
Also, neither the standard inner product of \eqref{typical-IP} nor the alternative inner product of \eqref{projection-IP} involve traces, whose existence in general depends on the dimensionality of their associated sets and properties of the involved weights (see, e.g., \cite[Th.~9.15]{Kufner:1985}).
Of course, one might use function restriction operators instead in certain classes of functions (e.g., polynomials, smooth functions), which, however, might fail to be complete with respect to the desired norm.
\end{remark}

Given $n \in \natural_0$, we denote by $\mathcal{V}^{\alpha,1}_n$ the space of $\HH^1_\alpha$-orthogonal polynomials of degree $n$ with respect to the $(\alpha,1)$-inner product given in \eqref{main-IP}; that is,
\begin{equation}\label{V1}
\mathcal{V}^{\alpha,1}_n := \left\{ p \in \poly^d_n \mid (\forall\,q\in \poly^d_{n-1}) \quad \langle p, q \rangle_{\alpha,1} = 0 \right\}.
\end{equation}
In common with the orthogonal polynomial spaces with respect to any inner product that is well defined over polynomials, $\mathcal{V}^{\alpha,1}_n$ inherits the dimension of the space of $d$-variate homogeneous polynomials of degree $n$ \cite[Sec.~3.1]{DunklXu:2014}:
\begin{equation}\label{dim-SOP}
\dim\left(\mathcal{V}^{\alpha,1}_n\right) = \binom{n+d-1}{n}.
\end{equation}
Also, as $\nabla 1 = 0$, for all $u \in \HH^1_\alpha$, $\langle u, 1 \rangle_{\alpha,1} = \langle S^\alpha_0(u), 1 \rangle_\alpha$, whence
\begin{equation}\label{SOP-Lebesgue-orthogonality}
(\forall\,n\in\natural) \quad \mathcal{V}^{\alpha,1}_n \perp_\alpha 1.
\end{equation}
Thus, we have the following analogue of \eqref{low-degree-OP}:
For all $\alpha > -1$,
\begin{equation}\label{low-degree-SOP}
\mathcal{V}^{\alpha,1}_0 = \vspan(\{1\})
\quad\text{and}\quad
\mathcal{V}^{\alpha,1}_1 = \vspan\left((x \mapsto x_i)_{i=1}^d\right).
\end{equation}

On account of \eqref{diff-shift}, one might expect $\mathcal{V}^{\alpha,1}_n$ to be related to $\mathcal{V}^{\alpha-1}_n$.
The following result states that this is indeed the case, with equality, if $\alpha > 0$ and $n \neq 2$.

\begin{proposition}\label{pro:high-parameter}
Let $d \in \natural$, $\alpha > 0$ and $n \in \natural_0 \setminus \{2\}$.
Then, $\mathcal{V}^{\alpha,1}_n = \mathcal{V}^{\alpha-1}_n$.
\begin{proof}
If $n = 0$ or $n = 1$ this comes from \eqref{low-degree-OP} and \eqref{low-degree-SOP}.
Let us suppose now that $n \geq 3$ and let $p_n \in \mathcal{V}^{\alpha-1}_n$.
Then, $\langle p_n, 1 \rangle_\alpha = \langle p_n, x \mapsto (1-\norm{x}^2) \rangle_{\alpha-1} = 0$, so $S^\alpha(p_n) = 0$.
Hence, given $q \in \poly^d_{n-1}$, $\langle p_n, q \rangle_{\alpha,1}
= \langle \nabla p_n, \nabla q \rangle_\alpha \stackrel{\eqref{diff-shift}}{=} 0$.
This establishes that $\mathcal{V}^{\alpha-1}_n$ is a subspace of $\mathcal{V}^{\alpha,1}_n$.
As, per \eqref{dim-OP} and \eqref{dim-SOP}, $\dim(\mathcal{V}^{\alpha-1}_n) = \dim(\mathcal{V}^{\alpha,1}_n)$, we obtain the desired equality.
\end{proof}
\end{proposition}

After a fashion, the rest of this work can be seen as an effort to extend the consequences of the above result to the whole natural range for $\alpha$, namely $(-1, \infty)$.

If $\alpha > 0$, from \eqref{parameter-lowering}, $d^{\alpha-1}_j \mathcal{V}^\alpha_n \subseteq \mathcal{V}^{\alpha-1}_{n+1}$.
Combining this with \eqref{diff-shift}, it is immediate that $\partial_i d^{\alpha-1}_j \mathcal{V}^\alpha_n \subseteq \mathcal{V}^\alpha_n$.
Even though the inclusion $d^{\alpha-1}_j \mathcal{V}^\alpha_n \subseteq \mathcal{V}^{\alpha-1}_{n+1}$ cannot be extended to $\alpha > -1$, its combination with \eqref{diff-shift} can, as we prove below in \autoref{pro:mappedGradientOrthogonality}.
First, however, we need the following commutation relations, the first of which already appears in the proof of \cite[Lemma~1.8.3]{DaiXu:2013}.

\begin{proposition}\label{pro:commutators}
Let $d \in \natural$, $i, j, k, l \in \{1, \dotsc, d\}$ and $\alpha \in \Real$.
Then,
\begin{equation}\label{diff-Dij-commutator}
\partial_i \, D_{k,l} - D_{k,l} \, \partial_i = \delta_{i,k} \partial_l - \delta_{i,l} \partial_k
\end{equation}
and
\begin{equation}\label{diff-PLO-commutator}
\partial_i d^{\alpha-1}_j - d^\alpha_j \partial_i = 2 D_{i,j} + 2\alpha \delta_{i,j} I,
\end{equation}
where $I$ is the identity operator.
\begin{proof}
Equations \eqref{diff-Dij-commutator} and \eqref{diff-PLO-commutator} are direct consequences of the definition of the operators $D_{k,l}$ and $d^\alpha_j$ in \eqref{Dij} and \eqref{PLO}, respectively, and the rules of calculus.
\end{proof}
\end{proposition}

\begin{proposition}\label{pro:mappedGradientOrthogonality}
Let $d \in \natural$, $\alpha > -1$ and $n \in \natural_0$.
\begin{enumerate}
\item\label{it:MGO-low} Let $i, j \in \{1, \dotsc, d\}$.
Then,
\begin{equation*}
\partial_i d^{\alpha-1}_j \mathcal{V}^\alpha_n \subseteq \mathcal{V}^\alpha_n.
\end{equation*}
\item\label{it:DijMGO-low} Let $i, j, k, l \in \{1, \dotsc, d\}$.
Then,
\begin{equation*}
\partial_i D_{k,l} d^{\alpha-1}_j \mathcal{V}^\alpha_n \subseteq \mathcal{V}^\alpha_n.
\end{equation*}
\end{enumerate}
\begin{proof}
Let $p_n \in \mathcal{V}^\alpha_n$.
By \eqref{diff-PLO-commutator} in \autoref{pro:commutators}, $\partial_i d^{\alpha-1}_j p_n$ coincides with $d^\alpha_j \partial_i p_n + 2 D_{i,j} p_n + 2\alpha\delta_{i,j} p_n$.
The first term belongs to $\mathcal{V}^\alpha_n$ because of \eqref{diff-shift} and \eqref{parameter-lowering}, the second because of \eqref{Dij-no-shift} and the third for free.
Thus, we have proved \autoref{it:MGO-low}.
Similarly, by \eqref{diff-Dij-commutator} in \autoref{pro:commutators}, $\partial_i D_{k,l} d^{\alpha-1}_j p_n$ coincides with $D_{k,l} \partial_i d^{\alpha-1}_j p_n + \delta_{i,k} \partial_l d^{\alpha-1}_j p_n - \delta_{i,l} \partial_k d^{\alpha-1}_j p_n$.
Each of the resulting three terms belongs to $\mathcal{V}^\alpha_n$ because of \autoref{it:MGO-low} with the help of \eqref{Dij-no-shift} in the case of the first.
This accounts for \autoref{it:DijMGO-low}.
\end{proof}
\end{proposition}

Let $M^\alpha$ denote the second-order differential operator defined by
\begin{equation}\label{M}
\begin{split}
M^\alpha(u)(x) & := (1-\norm{x}^2)^{1-\alpha} \Delta\left( (1-\norm{x}^2)^{1+\alpha} u(x) \right)\\
& \stackrel{\eqref{PLO}}{=} \sum_{j=1}^d d^{\alpha-1}_j d^\alpha_j u(x).
\end{split}
\end{equation}
In \autoref{lem:SOP} it will be proved that this operator transforms certain Lebesgue orthogonal polynomials into Sobolev orthogonal polynomials.
But first, we need some preliminary propositions and the inclusion, valid for all $\alpha > -1$, $n \in \natural_0$ and $i \in \{1, \dotsc, d\}$,
\begin{equation}\label{diff-SH}
\partial_i \mathcal{H}^d_n \subseteq \mathcal{H}^d_{n-1},
\end{equation}
which is easily derived from the fact that the partial derivative $\partial_i$ and the Laplacian $\Delta$ commute.

\begin{proposition}\label{pro:mutualHarmonicMappedGradientOrthogonality}
Let $d \in \natural$ and $\alpha > -1$.
\begin{enumerate}
\item\label{it:proto-MHMGO} Let $n \in \natural_0$, $h_n \in \mathcal{H}^d_n$ and $p_n \in \poly^d_n$.
Then,
\begin{equation*}
\langle \nabla h_n, \nabla p_n \rangle_\alpha = n(2n+d+2\alpha) \, \langle h_n, p_n \rangle_\alpha.
\end{equation*}
\item\label{it:main-MHMGO} Let $n \in \natural_0$, $h_n \in \mathcal{H}^d_n$ and $p_{n-2} \in \poly^d_{n-2}$.
Then,
\begin{equation*}
\langle h_n, M^\alpha(p_{n-2}) \rangle_\alpha = \langle \nabla h_n, \nabla M^\alpha(p_{n-2}) \rangle_\alpha = 0.
\end{equation*}
\end{enumerate}
\begin{proof}
We start by stating the fact that, for all $n \in \natural_0$ and $\mathcal{H}^d_n$,
\begin{equation}\label{SH-diagonalize-LB}
-\sum_{1\leq i<j\leq d} D_{i,j}^2 h_n = n(n+d-2) \, h_n.
\end{equation}
Indeed, as $h_n \in \mathcal{V}^0_n$ (cf.\ \eqref{SH-are-OP}), by \eqref{L}, \eqref{SL-0-strong} and \eqref{SL-0-strong-ew}, $n(n+d) \, h_n = -\Delta h_n + 2x \cdot \nabla h_n - \sum_{1 \leq i<j \leq d} D_{i,j}^2 h_n$; as $h_n$ is harmonic and homogeneous of degree $n$, $\Delta h_n = 0$ and $x \cdot \nabla h_n = n\,h_n$ \cite[Eq.~(12)]{DFLPP:2016}.
Alternatively, one might infer \eqref{SH-diagonalize-LB} by extending the corresponding result on the unit sphere, on which $\sum_{1\leq i<j\leq d} D_{i,j}^2$ becomes the Laplace--Beltrami operator \cite[Th.~1.4.5 and Th.~1.8.3]{DaiXu:2013}.

When $n = 0$, \autoref{it:proto-MHMGO} is obviously true.
So, let $n \geq 1$, $h_n \in \mathcal{H}^d_n$ and let us consider first the special case in which $p_n \in \mathcal{V}^\alpha_n$.
From \eqref{SH-are-OP}, $\mathcal{H}^d_n \subseteq \mathcal{V}^\alpha_n$.
By \eqref{diff-shift}, every component of both $\nabla h_n$ and $\nabla p_n$ belongs to $\mathcal{V}^{\alpha+1}_{n-1}$, so, by \cite[Prop.~3.4]{Figueroa:2017a},
\begin{equation}\label{wrongNorm}
\langle \nabla h_n, \nabla p_n \rangle_{\alpha+1} = \frac{\alpha+1}{n+d/2+\alpha} \, \langle \nabla h_n, \nabla p_n \rangle_\alpha.
\end{equation}
Also, by \autoref{pro:Dij-IBP} and \eqref{SH-diagonalize-LB},
\begin{equation}\label{weak-LB}
\sum_{1\leq i<j \leq d} \langle D_{i,j} h_n, D_{i,j} p_n \rangle_\alpha
= n(n+d-2) \langle h_n, p_n \rangle_\alpha.
\end{equation}
Combining \eqref{wrongNorm} and \eqref{weak-LB} with the weak form of the Sturm--Liouville problem that $h_n$ satisfies as a member of $\mathcal{V}^\alpha_n$ (cf.\ \eqref{L-self-adjoint} and \eqref{SL-0-weak}), we find that
\begin{equation*}
\lambda\sp{(\alpha)}_n \, \langle h_n, p_n \rangle_\alpha
= \frac{\alpha+1}{n+d/2+\alpha} \, \langle \nabla h_n, \nabla p_n \rangle_\alpha + n(n+d-2) \langle h_n, p_n \rangle_\alpha.
\end{equation*}
Rearranging the above equation and using the expression \eqref{SL-0-strong-ew} for $\lambda\sp{(\alpha)}_n$ accounts for \autoref{it:proto-MHMGO} in the special case in which $p_n \in \mathcal{V}^\alpha_n$.
From this special case the general $p_n \in \poly^d_n$ case follows upon the observation, arising from \eqref{SH-are-OP} and \eqref{diff-SH}, that $h_n$ is $\LL^2_\alpha$-orthogonal to $p_n - \proj^\alpha_n(p_n)$ and $\nabla h_n$ is $[\LL^2_\alpha]^d$-orthogonal to $\nabla(p_n - \proj^\alpha_n(p_n))$.

Let now $h_n \in \mathcal{H}^d_n$ and $p_{n-2} \in \poly^d_{n-2}$.
It is easy to check that
\begin{multline*}
\Delta \left[(1-\norm{x}^2)^{\alpha+1} \ p_{n-2}(x) \right]\\
= (1-\norm{x}^2)^{\alpha+1} \, \Delta p_{n-2}(x)
-4(\alpha+1) (1-\norm{x}^2)^\alpha \, x \cdot \nabla p_{n-2}(x)\\
+2(\alpha+1) \left[ 2\alpha (1-\norm{x}^2)^{\alpha-1} - (2\alpha+d) \, (1-\norm{x}^2)^\alpha \right] p_{n-2}(x).
\end{multline*}
Hence, using the first form of the operator $M^\alpha$ in \eqref{M},
\begin{multline*}
\langle h_n, M^\alpha(p_{n-2}) \rangle_\alpha
= \langle h_n, \Delta p_{n-2} \rangle_{\alpha+2}
- 4(\alpha+1) \langle h_n, x \cdot \nabla p_{n-2} \rangle_{\alpha+1}\\
+ 4\alpha(\alpha+1) \langle h_n, p_{n-2} \rangle_\alpha - 2(\alpha+1)(2\alpha+d) \langle h_n, p_{n-2} \rangle_{\alpha+1}.
\end{multline*}
As $h_n \in \mathcal{V}^\alpha_n \cap \mathcal{V}^{\alpha+1}_n \cap \mathcal{V}^{\alpha+2}_n$ (cf.\ \eqref{SH-are-OP}), $\langle h_n, M^\alpha(p_{n-2}) \rangle_\alpha$ vanishes and, by \autoref{it:proto-MHMGO}, so does $\langle \nabla h_n, \nabla p_n \rangle_{\alpha}$.
\end{proof}
\end{proposition}

\begin{proposition}\label{pro:mappedOrthogonality}
Let $d \in \natural$, $\alpha > -1$ and $n \in \natural_0$.
Then, $M^\alpha(\mathcal{V}^{\alpha+1}_n) \subseteq \mathcal{V}^\alpha_n \oplus_\alpha \mathcal{V}^\alpha_{n+2}$.
\begin{proof}
Let $p_n \in \mathcal{V}^{\alpha+1}_n$.
From the second form of $M^\alpha$ given in \eqref{M} and \eqref{PLO}, $M^\alpha(p_n) \in \poly^d_{n+2}$.
Now, for all $q \in \poly^d_{n-1}$, by integration by parts
\begin{equation*}
\langle M^\alpha(p_n), q \rangle_\alpha
= \int_{B^d} p_n(x) \, \overline{\Delta\left( (1-\norm{x}^2) \, q(x) \right)} \, (1-\norm{x}^2)^{\alpha+1} \dd x = 0,
\end{equation*}
where the absence of boundary terms is a consequence of the fact that $\alpha > -1$ and the vanishing of the latter integral comes about because the part under the conjugation bar is a polynomial of degree equal or less than $n-1$.
Thus, $M^\alpha(p_n) \perp_\alpha \poly^d_{n-1}$.
As the Laplacian operator and multiplication by centrally symmetric functions preserve the parity of a function, $M^\alpha(p_n)$ inherits the parity of $p_n$ given in \eqref{parity}, which, in turn, is the opposite of that of $\mathcal{V}^\alpha_{n+1}$, whence $M^\alpha(p_n) \perp_\alpha \mathcal{V}^\alpha_{n+1}$.
\end{proof}
\end{proposition}

We can now provide an orthogonal decomposition of the Sobolev orthogonal spaces $\mathcal{V}^{\alpha,1}_n$ closely related to that in \cite[Cor~2.4]{Xu:2008a} reproduced in \eqref{Xu-OD} (we discuss how closely in \autoref{ssc:MMO}).
Other related decompositions can be found in \cite[Rem.~2.1]{Xu:2008a} and \cite[Th.~4.1]{PinarXu:2009} and there is the hint of one in \cite[Eq.~(7.2)]{Xu:2017b}.

\begin{lemma}\label{lem:SOP}
Let $d \in \natural$ and $\alpha > -1$.
Then,
\begin{equation*}
\mathcal{V}^{\alpha,1}_n = \begin{cases}
\mathcal{V}^\alpha_n & \text{if } n \leq 2,\\
\mathcal{H}^d_n \oplus_{\alpha,1} \mathcal{M}^\alpha(\mathcal{V}^{\alpha+1}_{n-2}) & \text{if } n \geq 3.
\end{cases}
\end{equation*}
\begin{proof}
The cases $n = 0$ and $n = 1$ come directly from \eqref{low-degree-OP} and \eqref{low-degree-SOP}.
Now, given $p_2 \in \mathcal{V}^\alpha_2$ and $q \in \poly^d_1$, as $S^\alpha_0(p_2) = 0$, $\langle p_2, q \rangle_{\alpha,1} = \sum_{i=1}^d \langle \partial_i p_2, \partial_i q \rangle_\alpha$, which in turn vanishes because, by \eqref{diff-shift} and \eqref{low-degree-OP}, the $\partial_i p_2$ belong to $\mathcal{V}^{\alpha+1}_1 = \mathcal{V}^\alpha_1$.
Thus, $\mathcal{V}^\alpha_2 \subseteq \mathcal{V}^{\alpha,1}_2$ and as $\dim(\mathcal{V}^\alpha_2) = \dim(\mathcal{V}^{\alpha,1}_2)$ (cf.\ \eqref{dim-OP}, \eqref{dim-SOP}), the desired equality for the case $n = 2$ follows.

Let us suppose from now on that $n \geq 3$.
Let $h_n \in \mathcal{H}^d_n$.
From \eqref{SH-are-OP}, $S^\alpha_0(h_n) = 0$.
Hence, for all $q \in \poly^d_{n-1}$, $\langle h_n, q \rangle_{\alpha,1} = \sum_{i=1}^d \langle \partial_i h_n, \partial_i q \rangle_\alpha = 0$, the latter equality following from \eqref{SH-are-OP} and \eqref{diff-SH} on account of each of the $\partial_i q$ belonging to $\poly^d_{n-2}$.
Therefore, $h_n \in \mathcal{V}^{\alpha,1}_n$.

Let $p_{n-2} \in \mathcal{V}^{\alpha+1}_{n-2}$.
Then, $M^\alpha(p_{n-2}) \in \poly^d_n$ and, from \autoref{pro:mappedOrthogonality}, $S^\alpha_0(M^\alpha(p_n)) = 0$.
Also, from \eqref{parameter-lowering}, for every $j \in \{1, \dotsc, d\}$, $d^\alpha_j p_{n-2} \in \mathcal{V}^\alpha_{n-1}$.
Thus, for all $q \in \poly^d_{n-1}$,
\begin{equation*}
\langle M^\alpha(p_{n-2}), q \rangle_{\alpha,1}
= \sum_{i=1}^d \langle \partial_i M^\alpha(p_{n-2}), \partial_i q \rangle_\alpha
\stackrel{\eqref{M}}{=} \sum_{j=1}^d \sum_{i=1}^d \left\langle \partial_i d^{\alpha-1}_j d^\alpha_j p_{n-2}, \partial_i q \right\rangle_\alpha.
\end{equation*}
As, per \autoref{pro:mappedGradientOrthogonality}, each of the $\partial_i d^{\alpha-1}_j$ maps $\mathcal{V}^\alpha_{n-1}$ into itself and all the $\partial_i q$ belong to $\poly^d_{n-2}$, we find that $M^\alpha(p_{n-2}) \in \mathcal{V}^{\alpha,1}_n$.

Let $\{v_1, \dotsc, v_r\}$ be any basis of $\mathcal{V}^{\alpha+1}_{n-2}$ (by \eqref{dim-OP}, $r = \binom{n+d-3}{n-2}$) and let $a_1, \dotsc, a_r$ be scalars such that $\sum_{i=1}^r a_i \, \mathcal{M}^\alpha(v_i) = 0$.
Then, using the definition \eqref{M},
\begin{equation*}
(\forall\,x\in B^d) \quad \Delta\left( (1-\norm{x}^2)^{1+\alpha} \sum_{i=1}^r a_i \, v_i(x) \right) = 0.
\end{equation*}
Thus $x \mapsto (1-\norm{x}^2)^{1+\alpha} \sum_{i=1}^r a_i \, v_i(x)$ is harmonic in the unit ball and vanishes on the unit sphere.
As a consequence of the maximum principle for harmonic functions \cite[Th.~2.4]{GT}, it must be the null function and so must be $\sum_{i=1}^r a_i \, v_i$.
As the $v_i$ are mutually linearly independent, we conclude that $a_1 = \dotsb = a_r = 0$ and so $\{\mathcal{M}^\alpha(v_1), \dotsc, \mathcal{M}^\alpha(v_r)\}$ is a basis of $\mathcal{M}^\alpha(\mathcal{V}^{\alpha+1}_{n-2})$.
From \autoref{it:main-MHMGO} of \autoref{pro:mutualHarmonicMappedGradientOrthogonality} and the fact (cf.\ \eqref{SH-are-OP}) that $S^\alpha_0(\mathcal{H}^d_n) = \{0\}$, the vector space $M^\alpha(\mathcal{V}^{\alpha+1}_{n-2})$ is $(\alpha,1)$-orthogonal to $\mathcal{H}^d_n$.
Hence,
\begin{multline*}
\dim\left(\mathcal{H}^d_n \oplus \mathcal{M}^\alpha(\mathcal{V}^{\alpha+1}_{n-2})\right)
= \dim\left(\mathcal{H}^d_n\right) + \dim\left(\mathcal{M}^\alpha(\mathcal{V}^{\alpha+1}_{n-2})\right)\\
\stackrel{\eqref{dim-SH}}{=} \binom{n+d-1}{d-1} - \binom{n+d-3}{d-1} + \binom{n+d-3}{n-2}
= \binom{n+d-1}{n}
\stackrel{\eqref{dim-SOP}}{=} \dim\left(\mathcal{V}^{\alpha,1}_n\right).
\end{multline*}
As $\mathcal{H}^d_n \oplus \mathcal{M}^\alpha(\mathcal{V}^{\alpha+1}_{n-2}) \subseteq \mathcal{V}^{\alpha,1}_n$ with equality of dimensions, the proof is completed.
\end{proof}
\end{lemma}

\begin{remark}\label{rem:degree2}
By direct computation using the first form in \eqref{M}, for all $d \in \natural$ and $\alpha > -1$, there holds $\mathcal{M}^\alpha(1)(x) = 2(\alpha+1) \big((2\alpha+d)\norm{x}^2 - d\big)$.
With the help of known formulae \cite[Sec.~1.1 and Sec.~1.8]{AAR:1999} it is possible to evaluate $\langle M^\alpha(1), 1 \rangle_{\alpha,1} = \langle M^\alpha(1),1 \rangle_\alpha = -2 \, d \, \Gamma(1/2)^d \, \Gamma(\alpha+2) / \Gamma(d/2+\alpha+2) \neq 0$.
Therefore, for all $d \in \natural$ and $\alpha > -1$, the polynomial $M^\alpha(1)$ does not belong to $\mathcal{V}^{\alpha,1}_2$, so the latter orthogonal polynomial space does not fit into the decomposition pattern that \autoref{lem:SOP} sets for $\mathcal{V}^{\alpha,1}_n$, $n \geq 3$.
\end{remark}

We are now in position to provide a weak Sturm--Liouville problem that is satisfied by the $(\alpha,1)$-orthogonal polynomials.
As a preliminary step, we show that first-order partial derivatives of a Sobolev orthogonal polynomial in $\mathcal{V}^{\alpha,1}_n$ are Lebesgue orthogonal polynomials in $\mathcal{V}^\alpha_{n-1}$.

\begin{proposition}\label{pro:diff-SOP-is-OP}
Let $d \in \natural$, $\alpha > -1$ and $n \in \natural_0$.
Then, for all $k \in \{1, \dotsc, d\}$, $\partial_k \mathcal{V}^{\alpha,1}_n \subseteq \mathcal{V}^\alpha_{n-1}$.
\begin{proof}
If $n = 0$ or $n = 1$ this is immediate.
If $n = 2$, \autoref{lem:SOP} states that $p_n \in \mathcal{V}^\alpha_2$, so by \eqref{low-degree-OP} and \eqref{diff-shift}, $\partial_k p_n \in \mathcal{V}^\alpha_1$.
If $n \geq 3$, it transpires from \autoref{lem:SOP} that there exist $h_n \in \mathcal{H}^d_n$ and $r_{n-2} \in \mathcal{V}^{\alpha+1}_{n-2}$ such that $p_n = h_n + \mathcal{M}^\alpha(r_{n-2})$.
By \eqref{SH-are-OP} and \eqref{diff-SH}, $\partial_k h_n \in \mathcal{V}^\alpha_{n-1}$.
On the other hand, by the second form of $\mathcal{M}^\alpha$ given in \eqref{M}, $\partial_k \mathcal{M}^\alpha(r_{n-2}) = \sum_{l=1}^d \partial_k d^{\alpha-1}_l d^\alpha_l r_{n-2}$, so through \eqref{parameter-lowering} and \autoref{it:proto-MHMGO} of \autoref{pro:mappedGradientOrthogonality}, we infer that $\partial_k \mathcal{M}^\alpha(r_{n-2}) \in \mathcal{V}^\alpha_{n-1}$.
\end{proof}
\end{proposition}

\begin{thm}\label{thm:SL-1}
Let $d \in \natural$, $\alpha > -1$, $n \in \natural_0$ and $p_n \in \mathcal{V}^{\alpha,1}_n$.
Then, for all $q\in\CC^2(\overline{B^d})$,
\begin{multline}\label{SL-1}
B^{\alpha,1}(p_n,q) := \sum_{k=1}^d B^\alpha(\partial_k p_n, \partial_k q)\\
= \langle \nabla\nabla p_n, \nabla\nabla q \rangle_{\alpha+1} + \sum_{1 \leq i < j \leq d} \langle D_{i,j} \nabla p_n, D_{i,j} \nabla q \rangle_\alpha
= \lambda\sp{(\alpha,1)}_n \, \langle p_n, q \rangle_{\alpha,1},
\end{multline}
where
\begin{equation}\label{ew-1}
\lambda\sp{(\alpha,1)}_n
= \begin{cases}
0 & \text{if } n \leq 1,\\
\lambda\sp{(\alpha)}_{n-1} & \text{if } n \geq 2.
\end{cases}
\end{equation}
\begin{proof}
If $n \leq 1$, both $\nabla\nabla p_n$ and, for all admissible $i$ and $j$, $D_{i,j} \nabla p_n$ vanish and the desired result immediately follows.

From now on we suppose that $n \geq 2$.
From \autoref{pro:diff-SOP-is-OP}, $\partial_k p_n \in \mathcal{V}^\alpha_{n-1}$ for every $k \in \{1, \dotsc, d\}$.
As for every $q \in \CC^2(\overline{B^d})$, $\partial_k q \in \CC^1(\overline{B^d})$, we can substitute $n \leftarrow n-1$, $p_n \leftarrow \partial_k p_n$ and $q \leftarrow \partial_k q$ in the weak Sturm--Liouville problem \eqref{SL-0-weak}, sum up with respect to $k$ and obtain the desired result upon realizing that, on account of the Lebesgue orthogonality relation \eqref{SOP-Lebesgue-orthogonality}, $S^\alpha_0(p_n) = 0$ and hence $\langle \nabla p_n, \nabla q \rangle_\alpha = \langle p_n, q\rangle_{\alpha,1}$.
\end{proof}
\end{thm}

Given $\alpha > -1$, let us define a rank-1 perturbation $\tilde{\mathcal{L}}\sp{(\alpha)}$ of the second-order differential operator $\mathcal{L}\sp{(\alpha-1)}$ of \eqref{L} by
\begin{equation}\label{tilde-L}
(\forall\,q\in\HH^2_\alpha) \quad \tilde{\mathcal{L}}\sp{(\alpha)}(q) := \mathcal{L}\sp{(\alpha-1)}(q) + 2 \proj^\alpha_0(x \cdot \nabla q).
\end{equation}
Below we will prove that $\tilde{\mathcal{L}}\sp{(\alpha)}$ is $(\alpha,1)$-self-adjoint and that the strong Sturm--Liouville problem it defines is satisfied by the $(\alpha,1)$-orthogonal polynomials.
First, however, we need some relations connecting the $\mathcal{L}\sp{(\alpha)}$ and $\mathcal{L}\sp{(\alpha-1)}$ operators.

\begin{proposition}\label{pro:L-relations}
Let $d \in \natural$ and $\alpha \in \Real$.
Then,
\begin{equation}\label{L-difference}
\mathcal{L}\sp{(\alpha)} - \mathcal{L}\sp{(\alpha-1)} = 2 \, (x \cdot \nabla).
\end{equation}
and, for all $k \in \{1, \dotsc, d\}$,
\begin{equation}\label{L-derivative-commutator}
\mathcal{L}\sp{(\alpha)} \partial_k = \partial_k \mathcal{L}\sp{(\alpha-1)} - (d+2\alpha-1) \partial_k.
\end{equation}
\begin{proof}
The relation \eqref{L-difference} is a direct consequence of \eqref{PLO} and \eqref{L}.

By repeated application of \eqref{diff-Dij-commutator} of \autoref{pro:commutators}, for all $i, j$ satisfying $1 \leq i < j \leq d$ and $k \in \{1, \dotsc, d\}$,
\begin{equation*}
\begin{split}
D_{i,j}^2 \partial_k
& = D_{i,j} \left( \partial_k D_{i,j} - \delta_{k,i} \partial_j + \delta_{k,j} \partial_i \right)\\
& = \partial_k D_{i,j}^2 - \delta_{k,i} \partial_j D_{i,j} + \delta_{k,j} \partial_i D_{i,j} - \delta_{k,i} (\partial_j D_{i,j} + \partial_i) + \delta_{k,j} (\partial_i D_{i,j} - \partial_j)\\
& = \partial_k D_{i,j}^2 - 2 \delta_{k,i} \partial_j D_{i,j} + 2 \delta_{k,j} \partial_i D_{i,j} - \delta_{k,i} \partial_i - \delta_{k,j} \partial_j,
\end{split}
\end{equation*}
so
\begin{multline}\label{diff-LB-commutator}
-\sum_{1 \leq i < j \leq d} D_{i,j}^2 \partial_k\\
= -\partial_k \sum_{1 \leq i < j \leq d} D_{i,j}^2 + 2 \sum_{j=k+1}^d \partial_j D_{k,j} - 2 \sum_{i=1}^{k-1} \partial_i D_{i,k} + \sum_{j=k+1}^d \partial_k + \sum_{i=1}^{k-1} \partial_k\\
= -\partial_k \sum_{1 \leq i < j \leq d} D_{i,j}^2 + 2 \sum_{j=1}^d \partial_j D_{k,j} + (d-1) \partial_k.
\end{multline}
By \eqref{diff-Dij-commutator} and \eqref{diff-PLO-commutator} of \autoref{pro:commutators} and the fact that $\partial_j$ and $\partial_k$ commute, we find that, for all $j, k \in \{1, \dotsc, d\}$,
\begin{equation*}
\begin{split}
d^\alpha_j \partial_j \partial_k
& = \partial_k d^{\alpha-1}_j \partial_j - 2 D_{k,j} \partial_j - 2\alpha \delta_{k,j} \partial_j\\
& = \partial_k d^{\alpha-1}_j \partial_j - 2 (\partial_j D_{k,j} - \delta_{j,k} \partial_j + \partial_k) - 2\alpha \delta_{k,j} \partial_j,
\end{split}
\end{equation*}
whence
\begin{equation}\label{diff-wdivgrad-commutator}
\sum_{j=1}^d d^\alpha_j \partial_j \partial_k = \partial_k \sum_{j=1}^d d^{\alpha-1}_j \partial_j - 2 \sum_{j=1}^d \partial_j D_{k,j} - 2(d+\alpha-1) \partial_k.
\end{equation}
Summing \eqref{diff-LB-commutator} and \eqref{diff-wdivgrad-commutator} and using the first second form in \eqref{L} for both $\mathcal{L}\sp{(\alpha)}$ and $\mathcal{L}\sp{(\alpha-1)}$, we obtain \eqref{L-derivative-commutator}.
\end{proof}
\end{proposition}

\begin{lemma}\label{lem:tilde-L-self-adjoint}
Let $d \in \natural$ and $\alpha > -1$.
Then, the operator $\tilde{\mathcal{L}}\sp{(\alpha)}$ is self-adjoint in $\HH^1_\alpha$ with respect to the $(\alpha,1)$-inner product; indeed, for all $p \in \CC^3(\overline{B^d})$ and $q \in \CC^2(\overline{B^d})$,
\begin{equation*}
\left\langle \tilde{\mathcal{L}}\sp{(\alpha)}(p), q \right\rangle_{\alpha,1}
= B^{\alpha,1}(p,q) + (d+2\alpha-1) \langle \nabla p, \nabla q \rangle_\alpha.
\end{equation*}
\begin{proof}
Let us first notice that, as the gradient of any constant functions vanishes, $\nabla \tilde{\mathcal{L}}\sp{(\alpha)} = \nabla \mathcal{L}\sp{(\alpha-1)}$ (cf.\ \eqref{tilde-L}).
Using this fact, the equation \eqref{L-self-adjoint} encoding the self-adjointness of $\mathcal{L}\sp{(\alpha)}$ with respect to $\LL^2_\alpha$, the definition of the sesquilinear form $B^{\alpha,1}$ in \eqref{SL-1} and the relation \eqref{L-derivative-commutator} of \autoref{pro:L-relations}, we have
\begin{multline*}
\left\langle \nabla\tilde{\mathcal{L}}\sp{(\alpha)}(p), \nabla q \right\rangle_\alpha
= \left\langle \mathcal{L}\sp{(\alpha)}(\nabla p) + (d+2\alpha-1) \nabla p, \nabla q \right\rangle_\alpha\\
= B^{\alpha,1}(p,q) + (d+2\alpha-1) \langle \nabla p, \nabla q \rangle_\alpha.
\end{multline*}
Using the definition \eqref{tilde-L} of $\tilde{\mathcal{L}}\sp{(\alpha)}$, the self-adjointness of $\mathcal{L}\sp{(\alpha)}$ encoded in \eqref{L-self-adjoint} and the relation \eqref{L-difference} of \autoref{pro:L-relations}, we also have
\begin{equation*}
\proj^\alpha_0\left(\tilde{\mathcal{L}}\sp{(\alpha)}(p)\right)
= \proj^\alpha_0\left(\mathcal{L}\sp{(\alpha)}(p)\right)
= \frac{\langle \mathcal{L}\sp{(\alpha)}(p), 1 \rangle_\alpha}{\norm{1}_\alpha^2}
= \frac{\langle p, \mathcal{L}\sp{(\alpha)}(1) \rangle_\alpha}{\norm{1}_\alpha^2}
= 0,
\end{equation*}
and the desired result follows from the structure \eqref{main-IP} of the $(\alpha,1)$-inner product.
\end{proof}
\end{lemma}

\begin{thm}\label{thm:second-order-SL-1}
Let $d \in \natural$, $\alpha > -1$, $n \in \natural_0$ and $p_n \in \mathcal{V}^{\alpha,1}_n$.
Then,
\begin{equation*}
\tilde{\mathcal{L}}\sp{(\alpha)}(p_n) = \tilde\lambda\sp{(\alpha,1)}_n \, p_n,
\quad\text{where}\quad
\tilde\lambda\sp{(\alpha,1)}_n = n(n+d+2\alpha-2).
\end{equation*}
\begin{proof}
From \autoref{thm:SL-1} and \autoref{lem:tilde-L-self-adjoint}, for any polynomial $q$,
\begin{multline*}
\left\langle \tilde{\mathcal{L}}\sp{(\alpha)}(p_n), q \right\rangle_{\alpha,1}
= \lambda\sp{(\alpha,1)}_n \langle p_n, q \rangle_{\alpha,1} + (d+2\alpha-1) \langle \nabla p_n, \nabla q \rangle_\alpha\\
\stackrel{\text{\eqref{SOP-Lebesgue-orthogonality}}}{=} \left(\lambda\sp{(\alpha,1)}_n + \begin{cases} d+2\alpha-1 & \text{if $n \geq 1$},\\ 0 & \text{if $n = 0$} \end{cases}\right) \langle p_n, q \rangle_{\alpha,1}
\stackrel{\text{\eqref{SL-0-strong-ew},\eqref{ew-1}}}{=} \tilde\lambda\sp{(\alpha,1)}_n \langle p_n, q \rangle_{\alpha,1}.
\end{multline*}
As $\tilde{\mathcal{L}}\sp{(\alpha)}(p_n)$ itself is a polynomial (cf.\ \eqref{tilde-L}), the desired result follows.
\end{proof}
\end{thm}

\section{Approximation results}\label{sec:approximation}

Let us denote by $\proj^{\alpha,1}_n$ the $(\alpha,1)$-orthogonal projection from $\HH^1_\alpha$ onto $\mathcal{V}^{\alpha,1}_n$ \eqref{V1} and by $S^{\alpha,1}_n$ the $(\alpha,1)$-orthogonal projection from $\HH^1_\alpha$ onto $\poly^d_n$.
Clearly,
\begin{equation*}
\poly^d_n = \bigoplus_{k=0}^n \mathcal{V}^{\alpha,1}_k
\qquad\text{and}\qquad
S^{\alpha,1}_n = \sum_{k=0}^n \proj^{\alpha,1}_n.
\end{equation*}

\begin{proposition}\label{pro:poly-density}
Let $d \in \natural$ and $\alpha > -1$.
Then, polynomials are dense in $\HH^m_\alpha$.
\begin{proof}
By \autoref{pro:C-infty-density}, $\CC^\infty(\overline{B^d})$ is dense in $\HH^m_\alpha$; hence, so is the larger space $\CC^m(\overline{B^d})$.
Now, by \cite[Lem.~6.37]{GT}, every $\CC^m(\overline{B^d})$ function extends to a $\CC^m(\Real^d)$ function with compact support, so by \cite[Cor.~3]{EvardJafari:1994d} it can be approximated to within any positive distance by a polynomial in the norm $f \mapsto \sum_{\abs{\gamma}\leq m} \norm{\partial_\gamma f}_{\LL^\infty(B^d)}$, which is stronger than the norm of $\HH^m_\alpha$.
\end{proof}
\end{proposition}

As, by the above \autoref{pro:poly-density}, polynomials are dense in $\HH^1_\alpha$, there holds $\HH^1_\alpha = \bigoplus_{k=0}^\infty \mathcal{V}^{\alpha,1}_k$ so, by generic properties of Hilbert spaces \cite[Th.~I.4.13]{Conway:1990},
\begin{equation}\label{H1ExpansionParseval}
(\forall\,u\in\HH^1_\alpha)\quad u = \sum_{k=0}^\infty \proj^{\alpha,1}_k(u)
\quad\text{and}\quad
\norm{u}_{\alpha,1}^2 = \sum_{k=0}^\infty \norm[n]{\proj^{\alpha,1}_k(u)}_{\alpha,1}^2,
\end{equation}
the first series above converging in $\HH^1_\alpha$.

From its definition in \eqref{SL-1} within \autoref{thm:SL-1}, it is clear that the sesquilinear form $B^{\alpha,1}$ is still well defined on $\HH^2_\alpha$ and, with the help of the bound \eqref{B-bound} for $B^\alpha$, that it satisfies the bound
\begin{equation}\label{B1-bound}
(\forall\,u,v\in\HH^2_\alpha) \quad \abs{B^{\alpha,1}(u,v)} \leq \norm{\nabla\nabla u}_\alpha \norm{\nabla\nabla v}_\alpha;
\end{equation}
i.e., $B^{\alpha,1} \colon \HH^2_\alpha \times \HH^2_\alpha \to \Complex$ is a bounded sesquilinear form and hence, continuous.

We can now state one of our main approximation results, which provides bounds on the $\HH^1_\alpha$-best approximation error with respect to the $(\alpha,1)$-inner product of \eqref{main-IP} in terms of other best approximation errors.
This result lies at the foundation of the bound \eqref{simultaneous-H1} that will be attained later in \autoref{cor:old-style-bounds}.

\begin{thm}\label{thm:H1-approximation}
Let $d \in \natural$ and $\alpha > -1$.
Then,
\begin{multline}\label{proj-H1-grad-grad-BA}
(\forall\,u\in\HH^2_\alpha)\ (\forall\,N \geq 1)\\
\norm[n]{u - S^{\alpha,1}_N(u)}_{\alpha,1}
\leq (\lambda\sp{(\alpha,1)}_{N+1})^{-1/2} \inf_{p_N \in \poly^d_N} \norm{\nabla\nabla u - \nabla\nabla p_N}_\alpha
\end{multline}
and
\begin{multline}\label{proj-H1-SL-1-BA}
(\forall\,u\in\HH^3_\alpha)\ (\forall\,N \geq 2)\\
\norm[n]{u - S^{\alpha,1}_N(u)}_{\alpha,1} \leq (\tilde\lambda\sp{(\alpha,1)}_{N+1})^{-1} \norm{\tilde{\mathcal{L}}\sp{(\alpha)}(u) - S^{\alpha,1}_N(\tilde{\mathcal{L}}\sp{(\alpha)}(u))}_{\alpha,1}.
\end{multline}
\begin{proof}
We claim that, for all $v \in \HH^2_\alpha$ and $N \in \natural_0$,
\begin{equation}\label{minimizer-1}
S^{\alpha,1}_N(v) \in \operatorname*{arg\,min}_{q\in\poly^d_N} B^{\alpha,1}(v-q,v-q).
\end{equation}
and
\begin{equation}\label{Bessel-B1}
\sum_{n=0}^\infty \lambda\sp{(\alpha,1)}_n \norm{\proj^{\alpha,1}_n(v)}_{\alpha,1}^2 \leq B^{\alpha,1}(v,v).
\end{equation}
Indeed, \eqref{minimizer-1} and \eqref{Bessel-B1} are obtained by the same arguments that led to \eqref{minimizer} and \eqref{Bessel-B} in the proof of \autoref{thm:proj-L2-res-L2}, but with the $(\alpha,1)$-orthogonality taking the place of the $\LL^2$-orthogonality and substituting $B^\alpha \leftarrow B^{\alpha,1}$, \eqref{SL-0-weak} $\leftarrow$ \eqref{SL-1} and \eqref{B-bound} $\leftarrow$ \eqref{B1-bound}.

Thus, given $u \in \HH^2_\alpha$ and $N \geq 1$ (so that $\lambda\sp{(\alpha,1)}_n > 0$ for $n \geq N+1$; cf.\ \eqref{SL-0-strong-ew} and \eqref{ew-1}),
\begin{multline*}
\norm[n]{u - S^{\alpha,1}_N(u)}_{\alpha,1}^2 \leq \sup_{n \geq N+1} (\lambda\sp{(\alpha,1)}_n)^{-1} \sum_{n=N+1}^\infty \lambda\sp{(\alpha,1)}_n \norm{\proj^{\alpha,1}_n(u)}_{\alpha,1}^2\\
\stackrel{\text{\eqref{Bessel-B1}}}{\leq} (\lambda\sp{(\alpha,1)}_{N+1})^{-1} B^{\alpha,1}(u-S^{\alpha,1}_N(u), u-S^{\alpha,1}_N(u))\\
\stackrel{\text{\eqref{minimizer-1}}}{=} (\lambda\sp{(\alpha,1)}_{N+1})^{-1} \min_{p_N \in \poly^d_N} B^{\alpha,1}(u-p_N, u-p_N)
\stackrel{\text{\eqref{B1-bound}}}{\leq} (\lambda\sp{(\alpha,1)}_{N+1})^{-1} \inf_{p_N \in \poly^d_N} \norm{\nabla\nabla u - \nabla\nabla p_N}_\alpha^2,
\end{multline*}
whence \eqref{proj-H1-grad-grad-BA}.

Let $u \in \HH^3_\alpha$ and $N \geq 2$, so that $\tilde\lambda\sp{(\alpha,1)}_n > 0$ for $n \geq N+1$ (cf.\ \autoref{thm:second-order-SL-1}).
Continuing with the parallels with the proof of \autoref{thm:proj-L2-res-L2}, in the same way that \eqref{L-proj} is obtained there, here we obtain
\begin{equation}\label{L1-proj}
(\forall\,n\in\natural_0) \quad \proj^{\alpha,1}_n(\tilde{\mathcal{L}}\sp{(\alpha)} u) = \tilde\lambda\sp{(\alpha,1)}_n \proj^{\alpha,1}_n(u),
\end{equation}
where, instead of \eqref{SL-0-strong}, \eqref{SL-0-strong-ew}, \eqref{L-self-adjoint} and the continuity of $\mathcal{L}\sp{(\alpha)}$ as a map from $\HH^2_\alpha$ to $\LL^2_\alpha$ that are relevant for \eqref{L-proj}, we use \autoref{thm:second-order-SL-1}, \autoref{lem:tilde-L-self-adjoint} and the continuity of $\tilde{\mathcal{L}}\sp{(\alpha)}$ as a map from $\HH^3_\alpha$ to $\HH^1_\alpha$.
Then, the bound \eqref{proj-H1-SL-1-BA} follows from
\begin{multline*}
\norm{u - S^{\alpha,1}_N(u)}_{\alpha,1}^2
\leq \sup_{n \geq N+1} (\tilde\lambda\sp{(\alpha,1)}_n)^{-2} \sum_{n=N+1}^\infty \norm{\tilde\lambda\sp{(\alpha,1)}_n \proj^{\alpha,1}_n(u)}_{\alpha,1}^2\\
\stackrel{\text{\eqref{L1-proj}}}{=} (\tilde\lambda\sp{(\alpha,1)}_{N+1})^{-2} \sum_{n=N+1}^\infty \norm{\proj^{\alpha,1}_n(\tilde{\mathcal{L}}\sp{(\alpha)} u)}_{\alpha,1}^2\\
= (\tilde\lambda\sp{(\alpha,1)}_{N+1})^{-2} \norm{\tilde{\mathcal{L}}\sp{(\alpha)} u - S^{\alpha,1}_N(\tilde{\mathcal{L}}\sp{(\alpha)} u)}_{\alpha,1}^2.
\end{multline*}
\end{proof}
\end{thm}

\begin{remark}\label{rem:H1-approximation}
If $d+2\alpha > 0$, the bound \eqref{proj-H1-SL-1-BA} in \autoref{thm:H1-approximation} can start earlier at $N = 1$ and, if $d + 2\alpha > 1$, at $N = 0$.
\end{remark}

\autoref{pro:SOP-OP-id-shift} below echoes \cite[Prop.~3.1]{Figueroa:2017a} in establishing connections between orthogonal projection operators with respect to two different inner products.

\begin{proposition}\label{pro:SOP-OP-id-shift}
Let $d \in \natural$, $\alpha > -1$ and $n \in \natural_0$.
\begin{enumerate}
\item\label{it:SOP-OP-id-shift-basic} Let $q_n \in \mathcal{V}^{\alpha,1}_n$.
Then, $q_n = \proj^\alpha_{n-2}(q_n) + \proj^\alpha_n(q_n)$.
\item\label{it:SOP-OP-id-shift-meaty} Let $u \in \HH^1_\alpha$.
Then, $\proj^\alpha_n(u) = \proj^\alpha_n\left(\proj^{\alpha,1}_n(u)) + \proj^{\alpha,1}_{n+2}(u)\right)$.
\item\label{it:SOP-OP-id-shift-telescopable} Let $u \in \HH^1_\alpha$.
Then,
\begin{equation*}
\proj^\alpha_n(u) = \proj^{\alpha,1}_n(u) + \proj^\alpha_n \circ \proj^{\alpha,1}_{n+2}(u) - \proj^\alpha_{n-2} \circ \proj^{\alpha,1}_n(u).
\end{equation*}
\end{enumerate}
\begin{proof}
From \autoref{lem:SOP}, for $n \geq 3$, every $q_n \in \mathcal{V}^{\alpha,1}_n$ is of the form $q_n = h_n + r_n$, where $h_n \in \mathcal{H}^d_n \subseteq \mathcal{V}^\alpha_n$ (cf.\ \eqref{SH-are-OP}) and $r_n \in M^\alpha(\mathcal{V}^{\alpha+1}_{n-2}) \subseteq \mathcal{V}^\alpha_{n-2} \oplus_\alpha \mathcal{V}^\alpha_n$ (cf.\ \autoref{pro:mappedOrthogonality}).
The same \autoref{lem:SOP} states that, for $n \leq 2$, $\mathcal{V}^{\alpha,1}_n = \mathcal{V}^\alpha_n$.
Thus, for all $n \in \natural_0$, $\mathcal{V}^{\alpha,1}_n \subseteq \mathcal{V}^\alpha_{n-2} \oplus_\alpha \mathcal{V}^\alpha_n$, which implies \autoref{it:SOP-OP-id-shift-basic}.

Because of \eqref{H1ExpansionParseval}, the series $\sum_{k=0}^\infty \proj^{\alpha,1}_k(u)$ converges to $u$ in $\HH^1_\alpha$.
Because of the continuous embedding of $\HH^1_\alpha$ in $\LL^2_\alpha$, the same series converges to $u$ in $\LL^2_\alpha$ as well.
Then, \autoref{it:SOP-OP-id-shift-meaty} follows from
\begin{equation*}
(\forall\,p_n\in\mathcal{V}^\alpha_n) \quad \langle u, p_n \rangle_\alpha
= \sum_{k=0}^\infty \langle \proj^{\alpha,1}_k(u), p_n \rangle_\alpha
\stackrel{\text{\ref{it:SOP-OP-id-shift-basic}}}{=} \langle \proj^{\alpha,1}_n(u) + \proj^{\alpha,1}_{n+2}(u), p_n \rangle_\alpha.
\end{equation*}

Adding and subtracting $\proj^\alpha_{n-2}\left(\proj^{\alpha,1}_n(u)\right)$ to the right-hand side of \autoref{it:SOP-OP-id-shift-meaty} and using \autoref{it:SOP-OP-id-shift-basic} to combine $\proj^\alpha_n \circ \proj^{\alpha,1}_n(u)$ and $\proj^\alpha_{n-2} \circ \proj^{\alpha,1}_n(u)$ into $\proj^{\alpha,1}_n(u)$, we obtain \autoref{it:SOP-OP-id-shift-telescopable}.
\end{proof}
\end{proposition}

We now intend to prove another of our main approximation results, \autoref{thm:L2-approximation}, which bounds the $\LL^2$-error committed by the $(\alpha,1)$-best approximation in terms of the corresponding $(\alpha,1)$-best approximation error.
This result, in conjunction with \autoref{thm:H1-approximation}, lies at the foundation of the bound \eqref{simultaneous-L2} that will be attained later in \autoref{cor:old-style-bounds}.
First, we provide a bound on the difference between the $S^\alpha_N$ and the $S^{\alpha,1}_N$ projection operators in \autoref{lem:projectorDifference}.

\begin{lemma}\label{lem:projectorDifference}
Let $d \in \natural$ and $\alpha > -1$.
Then,
\begin{equation*}
(\forall\,u\in\HH^1_\alpha)\ (\forall\,N\geq 2) \quad \norm[n]{S^\alpha_N(u) - S^{\alpha,1}_N(u)}_\alpha \leq (\lambda\sp{(\alpha)}_{N-1})^{-1/2} \norm[n]{u - S^{\alpha,1}_N(u)}_{\alpha,1}.
\end{equation*}
\begin{proof}
For every $n \geq 3$ and $p_n \in \mathcal{V}^{\alpha,1}_n$,
\begin{multline*}
\norm{p_n}_\alpha^2
\stackrel{\text{Prop.\,\ref{pro:SOP-OP-id-shift}(\ref{it:SOP-OP-id-shift-basic})}}{=} \norm{\proj^\alpha_{n-2}(p_n)}_\alpha^2 + \norm{\proj^\alpha_n(p_n)}_\alpha^2\\
\stackrel{\text{\eqref{SL-0-weak}}}{=} \frac{1}{\lambda\sp{(\alpha)}_{n-2}} B^\alpha(\proj^\alpha_{n-2}(p_n),\proj^\alpha_{n-2}(p_n))  + \frac{1}{\lambda\sp{(\alpha)}_n} B^\alpha(\proj^\alpha_n(p_n),\proj^\alpha_n(p_n))\\
\leq \frac{1}{\lambda\sp{(\alpha)}_{n-2}} B^\alpha(p_n,p_n)
\stackrel{\text{\eqref{B-bound}}}{\leq} \frac{1}{\lambda\sp{(\alpha)}_{n-2}} \norm{\nabla p_n}_\alpha^2.
\end{multline*}
Hence, as long as $n \geq 2$,
\begin{equation}\label{residualTermsBound}
\norm{\proj^\alpha_{n-1} \circ \proj^{\alpha,1}_{n+1}(u)}_\alpha^2 \leq \norm{\proj^{\alpha,1}_{n+1}(u)}_\alpha^2 \leq \frac{1}{\lambda\sp{(\alpha)}_{n-1}} \norm{\proj^{\alpha,1}_{n+1}(u)}_{\alpha,1}^2.
\end{equation}

Now, for all $N \geq 0$,
\begin{multline*}
S^\alpha_N(u) - S^{\alpha,1}_N(u)
= \sum_{n=0}^N \left( \proj^\alpha_n(u) - \proj^{\alpha,1}_n(u) \right)\\
\stackrel{\text{Prop.~\ref{pro:SOP-OP-id-shift}(\ref{it:SOP-OP-id-shift-telescopable})}}{=} \sum_{n=0}^N \left( \proj^\alpha_n \circ \proj^{\alpha,1}_{n+2}(u) - \proj^\alpha_{n-2} \circ \proj^{\alpha,1}_n(u) \right)\\
= \proj^\alpha_N \circ \proj^{\alpha,1}_{N+2}(u) + \proj^\alpha_{N-1} \circ \proj^{\alpha,1}_{N+1}(u),
\end{multline*}
where, in the last equality, we have used the fact that the sum over $n$ telescopes and the fact that $\mathcal{V}^\alpha_{-2} = \mathcal{V}^\alpha_{-1} = \{0\}$.
Combining this with \eqref{residualTermsBound} we have that, for all $N \geq 2$,
\begin{multline*}
\norm{S^\alpha_N(u) - S^{\alpha,1}_N(u)}_\alpha^2
\leq \frac{1}{\lambda\sp{(\alpha)}_N} \norm{\proj^{\alpha,1}_{N+2}(u)}_{\alpha,1}^2 + \frac{1}{\lambda\sp{(\alpha)}_{N-1}} \norm{\proj^{\alpha,1}_{N+1}(u)}_{\alpha,1}^2\\
\leq \frac{1}{\lambda\sp{(\alpha)}_{N-1}} \norm{u - S^{\alpha,1}_N(u)}_{\alpha,1}^2.
\end{multline*}
The desired result follows upon taking square roots.
\end{proof}
\end{lemma}

\begin{thm}\label{thm:L2-approximation}
Let $d \in \natural$ and let $\alpha > -1$.
Then,
\begin{multline*}
(\forall\,u\in\HH^1_\alpha)\ (\forall\,N\geq 2)\\
\norm[n]{u - S^{\alpha,1}_N(u)}_\alpha \leq \left( (\lambda\sp{(\alpha)}_{N+1})^{-1/2} + (\lambda\sp{(\alpha)}_{N-1})^{-1/2} \right) \norm[n]{u - S^{\alpha,1}_N(u)}_{\alpha,1}
\end{multline*}
\begin{proof}
This follows from \eqref{proj-L2-grad-BA} of \autoref{thm:proj-L2-res-L2}, \autoref{lem:projectorDifference} and the obvious fact that $\inf_{p_{N} \in \poly^d_{N}} \norm{\nabla u - \nabla p_N}_\alpha \leq \norm[n]{\nabla u - \nabla S^{\alpha,1}_N(u)}_\alpha$.
\end{proof}
\end{thm}

\autoref{thm:H1-approximation} and \autoref{thm:L2-approximation} are analogous to \autoref{thm:proj-L2-res-L2} in the sense that they bound projection residuals in terms of powers of eigenvalues and best approximation errors.
We now turn them into bounds of the form \eqref{simultaneous}, which are given in terms powers of projection degrees and weighted Sobolev seminorms.

\begin{corollary}\label{cor:old-style-bounds}
Let $d \in \natural$ and $\alpha > -1$.
Then, for all integers $m \geq 1$, there exists $C = C(d,\alpha,m) > 0$ such that
\begin{equation}\label{proj-H1-res-H1-old-style}
(\forall\,u\in\HH^m_\alpha)\ (\forall\,N\geq m-1) \quad
\norm[n]{\nabla u - \nabla S^{\alpha,1}_N(u)}_\alpha
\leq C N^{1-m} \norm{\nabla^m u}_\alpha.
\end{equation}
and
\begin{equation}\label{proj-H1-res-L2-old-style}
(\forall\,u\in\HH^m_\alpha)\ (\forall\,N\geq \max(m-1,1)) \quad
\norm[n]{u - S^{\alpha,1}_N(u)}_\alpha
\leq C N^{-m} \norm{\nabla^m u}_\alpha.
\end{equation}
\begin{proof}
We begin by observing that, from its definition in \eqref{tilde-L}, it is straightforward to check that the operator $\tilde{\mathcal{L}}\sp{(\alpha)}$ belongs to $\mathcal{L}(\HH^{k+2}_\alpha,\HH^k_\alpha)$ for all $k \in \natural_0$.
On account of the norm equivalence given in \autoref{pro:projection-IP}, we are free to express this fact in the form
\begin{equation}\label{tilde-L-bounded}
(\forall\,k\in\natural_0)\ (\forall\,u\in\HH^{k+2}_\alpha) \quad \norm[n]{\tilde{\mathcal{L}}(u)}_{\alpha,k} \leq c_k \norm{u}_{\alpha,k+2},
\end{equation}
where $c_k > 0$ only depends on $\alpha$, $d$ and, of course, $k$.

Next, as the $S^{\alpha,1}_N$ are $(\alpha,1)$-orthogonal projections,
\begin{equation}\label{proto-1}
(\forall\,u\in\HH^1_\alpha)\ (\forall\,N\geq 0) \quad
\norm[n]{u - S^{\alpha,1}_N(u)}_{\alpha,1}
\leq \norm{u}_{\alpha,1}.
\end{equation}
From \eqref{proj-H1-grad-grad-BA} of \autoref{thm:H1-approximation} and bounding the infimum there by the result of taking the particular choice $p_N = 0$,
\begin{equation}\label{proto-2}
(\forall\,u\in\HH^2_\alpha)\ (\forall\,N\geq 1) \quad
\norm[n]{u - S^{\alpha,1}_N(u)}_{\alpha,1}
\leq (\lambda\sp{(\alpha,1)}_{N+1})^{-1/2} \norm{u}_{\alpha,2}.
\end{equation}
If $m$ is odd and $m \geq 3$, by applying \eqref{proj-H1-SL-1-BA} of \autoref{thm:H1-approximation} $\frac{m-1}{2}$ times, \eqref{proto-1} once and \eqref{tilde-L-bounded} $\frac{m-1}{2}$ times, we find that there exists $C_1 = C_1(d,\alpha,m) > 0$ such that
\begin{equation}\label{proto-odd}
(\forall\,u\in\HH^m_\alpha)\ (\forall\,N\geq 2) \quad
\norm[n]{u - S^{\alpha,1}_N(u)}_{\alpha,1}
\leq C_1 (\tilde\lambda\sp{(\alpha,1)}_{N+1})^{-\frac{m-1}{2}} \norm{u}_{\alpha,m}.
\end{equation}
If $m$ is even and $m \geq 4$, by applying \eqref{proj-H1-SL-1-BA} of \autoref{thm:H1-approximation} $\frac{m-2}{2}$ times, \eqref{proto-2} once and \eqref{tilde-L-bounded} $\frac{m-2}{2}$ times, we find that for these $m$ there also exists some $C_1 = C_1(d,\alpha,m) > 0$ such that
\begin{equation}\label{proto-even}
(\forall\,u\in\HH^m_\alpha)\ (\forall\,N\geq 2) \quad
\norm[n]{u - S^{\alpha,1}_N(u)}_{\alpha,1}
\leq C_1 (\tilde\lambda\sp{(\alpha,1)}_{N+1})^{-\frac{m-2}{2}} (\lambda\sp{(\alpha,1)}_{N+1})^{-1/2} \norm{u}_{\alpha,m}.
\end{equation}

From \eqref{SL-0-strong-ew}, \eqref{ew-1} of \autoref{thm:SL-1} and \autoref{thm:second-order-SL-1} and the restriction $\alpha > -1$, we find that $\lambda\sp{(\alpha,1)}_{N+1}$ is positive for all $N \geq 1$ and $\tilde\lambda\sp{(\alpha,1)}_{N+1}$ is positive for all $N \geq 2$.
As $\lim_{N\to\infty} \lambda\sp{(\alpha,1)}_{N+1} N^{-2} = \lim_{N\to\infty} \tilde \lambda\sp{(\alpha,1)}_{N+1} N^{-2} = 1$, it follows that there exists $C_2 = C_2(d,\alpha) > 0$ such that
\begin{equation}\label{ew-to-powers}
(\forall\,N\geq 1) \quad (\lambda\sp{(\alpha,1)}_{N+1})^{-1/2} \leq C_2 \, N^{-1}
\quad\text{and}\quad
(\forall\,N\geq 2) \quad (\tilde\lambda\sp{(\alpha,1)}_{N+1})^{-1/2} \leq C_2 \, N^{-1}.
\end{equation}
Combining \eqref{ew-to-powers} with \eqref{proto-1} for $m = 1$, \eqref{proto-2} for $m = 2$, \eqref{proto-odd} for odd $m \geq 3$ and \eqref{proto-even} for even $m \geq 4$ and exploiting the fact that each of the orthogonal projectors $S^{\alpha,1}_N$ leaves polynomials of degree $\leq N$ invariant, we find that, for all $m \geq 1$, there exist constants $C_3 = C_3(d,\alpha,m) > 0$ and $C_4 = C_4(d,\alpha,m) > 0$ such that
\begin{multline}\label{powers-and-seminorms}
(\forall\,u\in\HH^m_\alpha) \ (\forall\,N\geq m-1)\\
\norm[n]{u - S^{\alpha,1}_N(u)}_{\alpha,1}
= \norm[n]{u - S^\alpha_{m-1}(u) - S^{\alpha,1}_N\left(u - S^\alpha_{m-1}(u)\right)}_{\alpha,1}\\
\leq C_3 \, N^{1-m} \norm{u - S^\alpha_{m-1}(u)}_{\alpha,m}
\leq C_4 \, N^{1-m} \norm{\nabla^m u}_\alpha,
\end{multline}
from which \eqref{proj-H1-res-H1-old-style} follows.

From \eqref{SL-0-strong-ew}, in the same way we obtained \eqref{ew-to-powers}, we obtain that there exists $C_5 = C_5(d,\alpha) > 0$ such that
\begin{equation*}
(\forall\,N \geq 2) \quad (\lambda\sp{(\alpha)}_{N+1})^{-1/2} + (\lambda\sp{(\alpha)}_{N-1})^{-1/2} \leq C_5 \, N^{-1}.
\end{equation*}
Combining the above equation with \autoref{thm:L2-approximation} and \eqref{powers-and-seminorms} find that, for all $m \geq 1$, there exists a constant $C_6 = C_6(d,\alpha,m) > 0$ such that
\begin{equation*}
(\forall\,u\in\HH^m_\alpha) \ (\forall\,N\geq\max(m-1,2)) \quad
\norm[n]{u-S^{\alpha,1}_N(u)}_\alpha \leq C_6 \,  N^{-m} \norm{\nabla^m u}_\alpha.
\end{equation*}
This encompasses all instances of \eqref{proj-H1-res-L2-old-style} except for the two very special cases $(m,N) = (1,1)$ and $(m,N) = (2,1)$.
In the first special case, using the continuous embedding of $\HH^1_\alpha$ in $\LL^2_\alpha$ and the fact that $S^{\alpha,1}_1$ is an orthogonal projection, we find that there exists a constant $C_7 = C_7(d,\alpha) > 0$ such that
\begin{multline*}
(\forall\,u\in\HH^1_\alpha) \quad \norm[n]{u - S^{\alpha,1}_1(u)}_\alpha
= \norm[n]{u - S^\alpha_0(u) - S^{\alpha,1}_1\left(u - S^\alpha_0(u)\right)}_\alpha\\
\leq C_7 \, \norm[n]{u - S^\alpha_0(u)}_{\alpha,1}
= C_7 \, \norm{\nabla u}_\alpha.
\end{multline*}
For the second special case, in a very similar vein we find that there exists a constant $C_8 = C_8(d,\alpha) > 0$ such that
\begin{equation*}
(\forall\,u\in \HH^2_\alpha) \quad \norm[n]{u - S^{\alpha,1}_1(u)}_\alpha \leq C_8 \, \norm{\nabla\nabla u}_\alpha.
\end{equation*}
Hence, \eqref{proj-H1-res-L2-old-style} holds with $C = \max(C_6,C_7)$ if $m = 1$, $C = \max(C_6,C_8)$ if $m = 2$ and $C = C_6$ if $m \geq 3$.
\end{proof}
\end{corollary}

\section{Miscellaneous results and conclusion}\label{sec:miscellanea}

In this section we collect some results addressing certain features that appear recurrently in the literature cited in \autoref{ssc:related} but did not play an explicit role in the main body of the present work.

\subsection{A commutation relation involving differentiation and projection}\label{ssc:diff-proj-commutation}

Differentiation-projection commutation relations involving Sobolev orthogonalities have been obtained in \cite[Th.~4.5]{LiXu:2014}, \cite[Th.~5.6, Th.~7.5, etc.]{Xu:2017b} and \cite[Prop.~4.3 and Prop.~5.5]{MPPR:2023} and exploited in the pursuit of approximation results; see also \cite[Prop.~3.2(v)]{Figueroa:2017a} and \cite[Lem.~3.1]{PinarXu:2018} for variants involving Lebesgue orthogonalities only.
In this subsection we prove our own.

Given $d \in \natural$, $\alpha > -1$ and $i \in \{1, \dotsc, d\}$ and understanding $\HH^1_\alpha$ to be equipped with the $(\alpha,1)$-inner product, let $T^\alpha_i \colon \LL^2_\alpha \to \HH^1_\alpha$ denote the Hilbert space adjoint \cite[Def.~II.2.4]{Conway:1990} of the partial differentiation operator $\partial_i \colon \HH^1_\alpha \to \LL^2_\alpha$.
The operator $T^\alpha_i$ satisfies the relation
\begin{equation}\label{T}
(\forall\,f\in\LL^2_\alpha)\ (\forall\,v\in\HH^1_\alpha) \quad \langle v, T^\alpha_i(f) \rangle_{\alpha,1} = \langle \partial_i v, f \rangle_\alpha.
\end{equation}

\begin{proposition}\label{pro:undiscovered-OP-to-SOP}
Let $d \in \natural$, $\alpha > -1$, $i \in \{1, \dotsc, d\}$.
Then,
\begin{equation*}
(\forall\,n\in\natural_0) \quad T^\alpha_i(\mathcal{V}^\alpha_n) \subseteq \mathcal{V}^{\alpha,1}_{n+1}.
\end{equation*}
\begin{proof}
Let $p_n \in \mathcal{V}^\alpha_n$.
Given any $k \in \natural_0$ with $k \neq n+1$ and $q_k \in \mathcal{V}^{\alpha,1}_k$, we know from \autoref{pro:diff-SOP-is-OP} that $\partial_i q_k \in \mathcal{V}^\alpha_{k-1}$, whence, by \eqref{T}, $\langle q_k, T^\alpha_i(p_n) \rangle_{\alpha,1} = \langle \partial_i q_k, p_n \rangle_\alpha = 0$.
Hence, $\proj^{\alpha,1}_{n+1}(T^\alpha_i(p_n))$ is the only possibly non-vanishing term in the expansion \eqref{H1ExpansionParseval} of $T^\alpha_i(p_n)$ in terms of its projections onto Sobolev orthogonal polynomial spaces.
\end{proof}
\end{proposition}

\begin{lemma}\label{lem:diff-proj-commutation}
Let $d \in \natural$, $\alpha > -1$, $i \in \{1, \dotsc, d\}$ and $n \in \natural_0$.
Then,
\begin{equation*}
(\forall\,u\in\HH^1_\alpha) \quad
\partial_i \proj^{\alpha,1}_{n+1}(u) = \proj^\alpha_n(\partial_i u).
\end{equation*}
\begin{proof}
We know from \autoref{pro:diff-SOP-is-OP} that $\partial_i \proj^{\alpha,1}_{n+1}(u)$ belongs to $\mathcal{V}^\alpha_n$.
The desired result follows upon the observation that, for all $v_n \in \mathcal{V}^\alpha_n$,
\begin{equation*}
\langle \partial_i \proj^{\alpha,1}_{n+1}(u), v_n \rangle_\alpha
= \langle \proj^{\alpha,1}_{n+1}(u), T^\alpha_i(v_n) \rangle_{\alpha,1}
= \langle u, T^\alpha_i(v_n) \rangle_{\alpha,1}
= \langle \partial_i u, v_n \rangle_\alpha,
\end{equation*}
where the first and third equalities come from the operational definition of the map $T^\alpha_i$ in \autoref{pro:undiscovered-OP-to-SOP} and the second from the fact that $T^\alpha_i$ maps $\mathcal{V}^\alpha_n$ into $\mathcal{V}^{\alpha,1}_{n+1}$.
\end{proof}
\end{lemma}

\begin{remark}\label{rem:diff-proj-commutation}
We managed to prove some of our main results (e.g., \eqref{proj-H1-grad-grad-BA} of \autoref{thm:H1-approximation}) based on the commutation relation of \autoref{lem:diff-proj-commutation} without recourse to the characterization of Sobolev orthogonal polynomials as solutions of Sturm--Liouville problems, but could not do so for others (e.g., \eqref{proj-H1-SL-1-BA} of \autoref{thm:H1-approximation}).
\end{remark}

\subsection{Term-balancing constant}\label{ssc:balancing-constant}

Unlike many Sobolev inner products studied in the literature (including most of those cited in \autoref{ssc:related}), our $(\alpha,1)$-inner product of \eqref{main-IP} does not feature a term-balancing positive constant.
We will show here that, were we to introduce one, the resulting orthogonal polynomial spaces and corresponding orthogonal projection operators onto them would remain the same.
Indeed, given $d \in \natural$, $\alpha > -1$ and a term-balancing constant $\rho > 0$, let us define the $(\alpha,1,\rho)$-inner product on $\HH^1_\alpha$ through
\begin{equation*}
\langle u, v \rangle_{\alpha,1,\rho} := \rho \, \langle \nabla u, \nabla v \rangle_\alpha + \langle S^\alpha_0(u), S^\alpha_0(v) \rangle_\alpha.
\end{equation*}
Further, for all $n \in \natural_0$, we define the corresponding space of orthogonal polynomials of degree $n$ by
\begin{equation*}
\mathcal{V}^{\alpha,1,\rho}_n := \left\{ p \in \poly^d_n \mid (\forall\,q\in \poly^d_{n-1}) \quad \langle p, q \rangle_{\alpha,1,\rho} = 0 \right\}
\end{equation*}
and let $\proj^{\alpha,1,\rho}_n$ denote the $(\alpha,1,\rho)$-orthogonal projector from $\HH^1_\alpha$ onto $\mathcal{V}^{\alpha,1,\rho}_n$.

\begin{proposition}\label{pro:balancing-constant}
Let $d \in \natural$, $\alpha > -1$ and $\rho > 0$.
Then, for all $n \in \natural_0$, $\mathcal{V}^{\alpha,1,\rho} = \mathcal{V}^{\alpha,1}_n$ and $\proj^{\alpha,1,\rho}_n = \proj^{\alpha,1}_n$.
\begin{proof}
The equality $\mathcal{V}^{\alpha,1,\rho}_0 = \vspan(\{ 1 \}) = \mathcal{V}^{\alpha,1}_0$ is obvious.
So, let us suppose now that $n \geq 1$ and let $p_n \in \mathcal{V}^{\alpha,1,\rho}_n$.
For exactly the same reasons as in \eqref{SOP-Lebesgue-orthogonality}, there holds $\mathcal{V}^{\alpha,1,\rho}_n \perp_\alpha 1$.
In particular, $S^\alpha_0(p_n) = 0$, whence, for all $q \in \poly^d_{n-1}$,
\begin{equation*}
0 = \langle p_n, q \rangle_{\alpha,1,\rho}
= \rho \langle \nabla p_n, \nabla q \rangle_\alpha
= \rho \langle p_n, q \rangle_{\alpha,1}.
\end{equation*}
Thus, $\mathcal{V}^{\alpha,1,\rho}_n \subseteq \mathcal{V}^{\alpha,1}_n$.
As both the involved inner products are well-defined on $d$-variate polynomials, $\dim(\mathcal{V}^{\alpha,1,\rho}_n) = \dim(\mathcal{V}^{\alpha,1}_n)$, so actually $\mathcal{V}^{\alpha,1,\rho}_n = \mathcal{V}^{\alpha,1}_n$.

Let $u \in \HH^1_\alpha$ and $n \in \natural_0$.
Then, $u^{(\rho)}_n := \proj^{\alpha,1,\rho}_n(u)$ is characterized by the pair of statements (\textsc{i}) $u^{(\rho)}_n \in \mathcal{V}^{\alpha,1,\rho}_n$ and (\textsc{ii}) $\langle u - u^{(\rho)}_n, q_n \rangle_{\alpha,1,\rho} = 0$ for all $q_n \in \mathcal{V}^{\alpha,1,\rho}_n$.
From the above equality of orthogonal polynomial spaces, it is immediate that (\textsc{i'}) $u^{(\rho)}_n \in \mathcal{V}^{\alpha,1}_n$.
Further, using (\textsc{i}) and the same equality of orthogonal polynomial spaces, we have that, for all $q_n \in \mathcal{V}^{\alpha,1}_n$,
\begin{multline*}
0 = \langle u - u^{(\rho)}_n, q_n \rangle_{\alpha,1,\rho} = \rho \langle \nabla u - \nabla u^{(\rho)}_n, \nabla q_n \rangle_\alpha + \langle S^\alpha_0(u - u^{(\rho)}_n), S^\alpha_0(q_n) \rangle_\alpha\\
= \langle u - u^{(\rho)}_n, q_n \rangle_{\alpha,1} \times
\begin{cases}
1 & \text{if } n = 0,\\
\rho & \text{if } n \geq 1.
\end{cases}
\end{multline*}
because $\nabla q_n = 0$ if $n = 0$ and because, as mentioned above, $q_n \in \mathcal{V}^{\alpha,1}_n = \mathcal{V}^{\alpha,1,\rho}_n \perp_\alpha 1$ if $n \geq 1$.
Hence, we have (\textsc{ii'}) $\langle u - u^{(\rho)}_n, q_n \rangle_{\alpha,1} = 0$ for all $q_n \in \mathcal{V}^{\alpha,1}_n$, which, in conjunction with (\textsc{i'}), ensures that $u^{(\rho)}_n = \proj^{\alpha,1}_n(u)$.
\end{proof}
\end{proposition}

\subsection{Multiplication operator in the unweighted case}\label{ssc:MMO}

In the unweighted ($\alpha = 0$) case, the operator $M^0$ of \eqref{M}, which, by \autoref{lem:SOP}, maps the Lebesgue orthogonal polynomial space $\mathcal{V}^1_{n-2}$ onto a subspace of the Sobolev orthogonal polynomial space $\mathcal{V}^{0,1}_n$, can be substituted with the multiplication-by-$(1-\norm{x}^2)$ operator.
As a consequence, for $n \neq 2$, the orthogonal decomposition of $\mathcal{V}^{0,1}_n$ given in \autoref{lem:SOP} is exactly the decomposition of the orthogonal polynomials with respect to the inner product \eqref{LiXu-IP} given in \cite[Cor.~2.4]{Xu:2008a} and reproduced in \eqref{Xu-OD}.

\begin{proposition}\label{pro:MMO}
Let $d \in \natural$.
\begin{enumerate}
\item\label{it:MMO} For all $n\in\natural_0$, $\mathcal{M}^0(\mathcal{V}^1_n) = (1-\norm{x}^2) \, \mathcal{V}^1_n$.
\item\label{it:same-SOP} For all $n\in\natural_0\setminus\{2\}$,
\begin{equation}\label{same-SOP}
\begin{split}
\mathcal{V}^{0,1}_n & = \mathcal{H}^d_n \oplus (1-\norm{x}^2) \, \mathcal{V}^1_{n-2}\\
& = \left\{ p \in \poly^d_n  \mid (\forall\,q\in\poly^d_{n-1})\ p \perp q \text{ with respect to \eqref{LiXu-IP}} \right\}.
\end{split}
\end{equation}
\end{enumerate}
\begin{proof}
For all $u \in C^2(B^d)$,
\begin{multline*}
\Delta\left( (1-\norm{x}^2) \, u(x) \right)
= \div\left( (1-\norm{x}^2) \nabla u(x) - 2 u(x) x \right)\\
\stackrel{\text{\eqref{L}}}{=} -\mathcal{L}^{(0)}(u)(x) - \sum_{1\leq i<j\leq d} (x_i \partial_j - x_j \partial_i)^2 u(x) - 2 \, x \cdot \nabla u(x) - 2 d \, u(x)\\
\stackrel{\text{\eqref{L-difference}}}{=} -\mathcal{L}^{(1)}(u)(x) - \sum_{1\leq i<j\leq d} (x_i \partial_j - x_j \partial_i)^2 u(x) - 2 d \, u(x).
\end{multline*}
This shows, via \eqref{SL-0-strong} and \eqref{Dij-no-shift}, that $\Delta\left( (1-\norm{x}^2) \, \mathcal{V}^1_n \right) \subseteq \mathcal{V}^1_n$.
If $u \in \mathcal{V}^1_n$ is such that $\Delta\left( (1-\norm{x}^2) \, u(x) \right) = 0$, then $u$ is the null function as a consequence of the maximum principle for harmonic functions \cite[Th.~2.4]{GT}.
Then, from the rank-nullity theorem of Linear Algebra \cite[Th.~50/1]{Halmos:1958}, $\Delta\left( (1-\norm{x}^2) \, \mathcal{V}^1_n \right) = \mathcal{V}^1_n$.
Thus, \autoref{it:MMO} follows from the first form of $M^0$ given in \eqref{M}.

The first equality in \eqref{same-SOP} is simply \autoref{lem:SOP} in combination with \autoref{it:MMO} (if $n \geq 3$) and with \eqref{low-degree-SOP} (if $n \leq 1$).
The second equality in \eqref{same-SOP} is the decomposition \eqref{Xu-OD}, trivially extended to the $n = 0$ case, and this accounts for \autoref{it:same-SOP}.
\end{proof}
\end{proposition}

\begin{remark}\label{rem:MMO}
If $n = 2$, the second equality of \eqref{same-SOP} in \autoref{it:same-SOP} of \autoref{pro:MMO} is true by \eqref{Xu-OD}, but the first equality is false on account of \autoref{it:MMO} and \autoref{rem:degree2}.
\end{remark}

\subsection{Conclusion}\label{ssc:conclusion}

We have proved our desired simultaneous approximation bounds \eqref{simultaneous} in \autoref{cor:old-style-bounds}.
To do so, we have delved into the rich structure of orthogonal polynomials with respect to the $(\alpha,1)$-inner product of \eqref{projection-IP-1} and \eqref{main-IP} arguing, as in \cite{Figueroa:2017a}, mostly in terms of orthogonal polynomial spaces instead of particular bases of them.
We expect that these results and techniques will generalize readily to other domains and families of symmetry-respecting weights.

\bibliographystyle{amsplain}
\bibliography{sop-refs}

\end{document}